\documentclass[11pt,a4paper]{article}
\usepackage{amsmath}
\usepackage[latin1]{inputenc}
\usepackage{amsfonts}
\usepackage{amssymb}
\usepackage{graphicx}
\usepackage{epsfig}
\usepackage[onehalfspacing]{setspace}
\newtheorem{theorem}{Theorem}[section]
\newtheorem{corollary}{Corollary}[section]

\newtheorem{lemma}{Lemma}[section]
\newtheorem{proposition}{Proposition}[section]
\newtheorem{remark}{Remark}[section]
\newenvironment{proof}[1][Proof]{\textbf{#1.} }{\ \rule{0.5em}{0.5em}}
\begin{document}

\title{CLT in Functional Linear Regression Models}
\author{Herv\'e Cardot$^{1,3}$, Andr\'e Mas$^{2}$, Pascal Sarda$^{3,4}$\\
{\small (1) INRA Toulouse, Biom\'etrie et Intelligence Artificielle,} \\
{\small31326 Castanet-Tolosan, France} \\
{\small (2) Institut de modélisation mathématique de Montpellier 2, Universit\'e Montpellier II,}\\
{\small Place Eugène Bataillon, 34095 Montpellier Cedex 5, France} \\
{\small (3) Laboratoire de Statistique et Probabilit\'es, Universit\'e Paul Sabatier} \\
{\small 118, route de Narbonne, 31062 Toulouse Cedex 4, France}\\
{\small (4) GRIMM, EA 2254, Université Toulouse-le-Mirail}\\
{\small 5, Allées Antonio-Machado, 31058 Toulouse Cedex 1, France}}
\maketitle

\begin{abstract}
We propose in this work to derive a CLT in the functional linear regression model.
The main difficulty is due to the fact that estimation of the functional parameter leads to a kind of ill-posed inverse problem. We consider estimators that belong to a large class of regularizing methods and we first show that, contrary to the multivariate case, it is not possible to state a CLT in the topology of the considered functional space.  
However, we show that we can get a CLT for the weak topology under mild hypotheses and in particular without assuming any  strong assumptions on the decay of the eigenvalues of the covariance operator. Rates of convergence depend on the smoothness of the functional coefficient and on the point in which the prediction is made.
\end{abstract}



\section{Introduction}


For several years, there has been a considerable interest in {\it Functional Data Analysis}.  Indeed, a consequence of advances in  technology is the collection of many data sets on dense grids (image, satellite, medicine, ...) adding in some sense more and more information. The question is then: can we do something specific with this new information~? It is the merit of the books by Ramsay and Silverman (1997, 2002) to have prepared the ground for answers to this question. They, and other authors after them, have shown the practical benefits of using {\it ad hoc} statistical methods for these data sets: the point of view is clearly to take into account the functional nature of the data. This means that one considers the data as objects belonging to functional spaces with infinite dimension and one has to use adapted probabilistic and functional analysis tools to derive properties of estimators in such  a context.

This emulates a need for developing theoretical/practical aspects on the ground of functional data analysis. It is the aim of this paper to contribute to this kind of development. The framework of our study is itself an important part of functional data problems. 
We are interested in the  properties of the linear regression model in this functional framework, that is to say performing the regression of a real random variable on a functional variable. The two main motivations of this work are to study rigorously the asymptotic distribution of the estimator of the regression function and from a statistical point of view to deduce asymptotic confidence intervals for prediction based on functional linear regression.

This kind of model is not  new and has many potential applications such as Chemometrics as it can be noticed in the paper by Frank and Friedman (1993). Whereas chemometricians have mainly used adaptations of statistical multivariate methods, {\it functional} procedures have gained in popularity more recently as said above. 
For instance Hastie and Mallows (1993) have raised, in the discussion of the paper by Frank and Friedman (1993), the question of functional alternative methods. Thus, for this case of estimating a regression, two main approaches have been considered, (1) estimating the functional linear regression which is a ``continuous'' version of linear regression when the covariate is a vector of scalars and was first introduced in Ramsay and Dalzell (1991) and (2) proposing a complete nonparametric point of view introduced by Ferraty and Vieu (2002). 
We consider the former approach hereafter: see section 2, for the definition of the functional linear regression. Contrary to the  multivariate linear regression where the vector of parameters is identifiable provided that the covariance matrix is non-singular, identifiability of the functional coefficient is not ensured unless a sufficient and necessary condition is satisfied (see section 2).

Different estimators for the functional parameter have been considered in the functional linear regression model: see for instance Goutis (1998), Cuevas {\it et al.} (2002) and Cardot {\it et al.} (1999, 2003). 
Upper bounds for the $L^2$ rates of convergence have been found and these results show that the transition from the finite dimension to the infinite dimension leads to degenerate rates of convergence. 
As a matter of fact, estimating the functional parameter appearing in the functional linear regression can be viewed as an ill-conditioned inverse problem since the estimation procedure  relies on the inversion of the covariance operator which is compact (Cardot \textit{et al.} 1999, Bosq 2000, He \textit{et al.} 2003) and it is well known that poor rates of convergence appear for this kind of problems. However the problem of approximating inverses of covariance operators or of selfadjoint compact operators is not new. It is adressed in Nashed and Wahba (1974), Arsenin and Tikhonov (1977), Groetsch (1993) among many others. The main point is always to regularize a matrix $M$ (respectively an operator $S$) which is invertible but ``not by much'' (respectively unbounded). This property implies that for any vector $x$, $Mx$ (respectively $Sx$) may have large variations even when $x$ does not vary much. Numerous procedures were proposed. Such procedures appear especially in image analysis or deconvolution or in specific M-estimation problems for instance.

The Central Limit Theorem for i.i.d. Hilbert valued random variables play a central role in deriving the main results of this paper. The monograph by Araujo and Giné (1980) or Chapter 5 in the book by Ledous and Talagrand (1991) deal with this crucial theorem of probability theory and provide deep studies about the CLT on infinite dimensional spaces. For the non independent case we also mention two recent works from Dedecker and Merlevède (2002) and Merlevède (2003).

In section 3, we consider a class of regularization methods for inverting the covariance operator that leads to a quite general class of estimators with the aim of investigating CLT for prediction and as a by-product producing confidence sets for prediction. 
Section 4 is devoted to the asymptotic behavior of these estimators  relaxing as much as possible the set of assumptions (moment assumptions, assumptions on the spectrum of $\Gamma$) and considering a large class of regularizing methods for inverting the empirical covariance operator.
We first derive an important result which shows that it is not possible to state a CLT for the functional coefficient with respect to the norm topology of the functional space. Nevertheless, we show that it is possible to  get a CLT if we consider the behavior of the predictor with respect to the weak topology, that is to say for point-wise prediction. 
We show that the results depend on the nature of the predictor and fixed or random design lead to different CLT. Whereas when the predictor is random it is not possible to reach a parametric rate of convergence, this rate can be obtained depending on the value and the smoothness properties of the fixed predictor: we obtain a parametric rate for pointwise convergence at $x$ wherver $x$ belongs to the reproducing kernel Hilbert space associated to the covariance operator. The proofs depend heavily on perturbation theory for linear operators to get, as accurate as possible, approximations of the eigenelements of the empirical covariance operators. Similar methods based on functional calculus have been used for deriving asymptotic properties of the functional  principal components analysis by Dauxois \textit{et al.} (1982),  Kneip (1995) or Kneip and Utikal (2001). 
 Section 5 proposes a brief discussion about possible extensions and statistical applications of these results. 
Finally section 6 is devoted to the proofs.


\section{Functional linear regression}


We consider a sample $(X_{i},Y_{i}),i=1,\ldots,n$ of independent and identically distributed random variables drawn from a pair $(X,Y)$. The variables $X$ and $Y$ are defined on the same probability space and $Y$ (the response) is valued in $\mathbb{R}$. It is usual to define the functional variable $X$ (the predictor) as a random variable taking values in a general real separable Hilbert space $H$ with an inner product denoted in the following by $\langle.,.\rangle$ and an associated norm denoted by $\Vert.\Vert$. As a matter of fact $H$ may be the Sobolev space $W^{m,2}\left(  \mathcal{C}\right)  $ of functions defined on some compact interval $\mathcal{C}$ of $\mathbb{R}$ having $m$ square integrable derivatives, $m$ being a positive integer. In that case the inner product $\langle.,.\rangle$ is the usual inner product on this space \textit{i.e.}
$$
\langle f,g\rangle\ =\ \sum_{p=0}^{m}\int_{\mathcal{C}}f^{\left(  p\right)}(x)g^{\left(  p\right)  }(x)dx,\qquad f,g\in H.
$$
Note that this special case is particularly interesting for modelling situations where we have functional data as shown by the numerous applications given in Ramsay and Silverman (1997, 2002). Although we develop below theory for general Hilbertian random variables, we keep in mind this special situation and then use the word \textit{functional} variable to qualify $X$. 

In the following we assume that $I\!\!EY^2<+\infty$ and that $X$ is a $H$-valued random variable such that
$$
I\!\!E(\Vert X\Vert^{4})\ <\ +\infty.\leqno(H.1)
$$
Then $X$ is of second order and one can define the expectation of $X$, namely $I\!\!E(X)$, that we suppose in order to simplify the notations to be the null element of $H$ ($I\!\!E(X)=0$). Moreover the covariance operator of $X$ is defined as the linear operator $\Gamma$ defined on $H$ as follows
$$
\Gamma h\ =\ I\!\!E(X\otimes X(h))\ =\ I\!\!E(\langle h,X\rangle X),\qquad
h\in H.
$$
It is known that $\Gamma$ is a self-adjoint, positive and nuclear operator hence it is Hilbert-Schmidt and hence compact (Dauxois {\it et al.}, 1982). We denote by $(\lambda_{j})_{j}$ the sorted sequence of non null distinct eigenvalues of $\Gamma$, $\lambda_1>\lambda_2>\ldots>0$, and $(e_{j})_{j}$ a sequence of orthonormal associated eigenvectors. We assume that the multiplicity of each $\lambda_{j}$ is one (remind that since $\Gamma$ is compact the multiplicity of each $\lambda_j\neq 0$ is finite). We could consider the more general case of multiple eigenvalues without affecting our forthcoming results but the price would be more complicated proofs and a poor gain with respect to the main objectives of the paper. Let us also define the cross-covariance operator of $X$ and $Y$ as the functional $\Delta$ defined on $H$ by
$$
\Delta h\ =\ I\!\!E(X\otimes Y(h))\ =\ I\!\!E(\langle h,X\rangle Y),\qquad h\in H.
$$

Now, we aim at considering the \textit{functional linear regression} of the variable $Y$ on $X$. This means that we are seeking the solution $\rho\in H$ of the following minimization problem
\begin{equation}\label{infbeta}
\inf_{\beta\in H}I\!\!E\left(  |Y-\langle\beta,X\rangle|^{2}\right).
\end{equation}
When a solution $\rho$ exists and is uniquely determined, we can write
\begin{equation}
Y\ =\ \langle\rho,X\rangle\,+\,\varepsilon, \label{model}
\end{equation}
where $\varepsilon$ is a centered real random variable with variance $\sigma_{\varepsilon}^{2}$ such that $E(\varepsilon X)=0$. It is quite easy to show that it is equivalent that $\rho$ satisfies equation (\ref{model}) and that it satisfies the following moment equation
$$
\label{normalequation}\Delta\ =\ \Gamma\rho.
$$
However, when the dimension of $H$ is infinite, existence and uniqueness of $\rho$ is not ensured
since a bounded inverse of $\Gamma$ does not exist:
we need an additional condition to get existence and uniqueness of $\rho$, namely

\vspace{.5cm}

\noindent\textbf{Condition }$\mathcal{U}$. The variables $X$ and $Y$ satisfy
$$
\sum_{j}\frac{\langle I\!\!E(XY),e_{j}\rangle^{2}}{\lambda_{j}^{2}}\ <\ +\infty.
$$

\vspace{.5cm}

Under condition $\mathcal{U}$, Cardot \textit{et al.} (2003) show that a unique solution to equation (\ref{model}) exists in $((Ker(\Gamma))^{\perp}$ and that this solution is of the form
$$
\rho\ =\ \sum_{j}\frac{\langle I\!\!E(XY),e_{j}\rangle}{\lambda_{j}}e_{j}.
$$
Then, identifiability is true only in $((Ker(\Gamma))^{\perp}$ or in other words the set of solution of (\ref{model}) is of the form $\rho+Ker(\Gamma)$. Again, to simplify further developments we assume from now on that the
following condition is satisfied
$$
Ker(\Gamma)\ =\ \{0\}.\leqno(H.2)
$$
Finally, we assume from now on that the first and second moment of $\varepsilon$ given $X$ are respectively equal to $I\!\!E(\varepsilon|X)=0$ and $I\!\!E(\varepsilon^2|X)=\sigma^2_\varepsilon$.


\section{Inverse problem and regularization procedure}


Once we get identifiability through condition $\mathcal{U}$, we turn to the problem of estimating the ``functional'' parameter $\rho$ from the sample $(X_{i},Y_{i}),\ i=1,\ldots,n$. The first step is to define the empirical versions of $\Gamma$ and $\Delta$ which are
$$
\Gamma_{n}\ =\ \frac{1}{n}\sum_{i=1}^{n}X_{i}\otimes X_{i},\qquad
\Delta_{n}\ =\ \frac{1}{n}\sum_{i=1}^{n}X_{i}\otimes Y_{i}.
$$
We have
$$
\Delta_{n}\ =\ \Gamma_{n}\rho\,+\,U_{n},
$$
where $U_{n}=n^{-1}\sum_{i=1}^{n}X_{i}\otimes\varepsilon_{i}$ and taking the expectation we get
$$
I\!\!E(\Delta_{n})\ =\ \Delta\ =\ \Gamma\rho.
$$
As shown in the previous section, inversion of $\Gamma$ can be viewed as a kind of ill-conditioned inverse problem (unlike in usual ill-conditioned inverse problems the operator $\Gamma$ is unknown). Also, the inverse of $\Gamma_{n}$ does not exist because $\Gamma_{n}$ is almost surely a finite rank operator. As usually for ill-conditioned inverse problem we need regularization and our aim is now to propose a general and unified method to get a sequence of continuous estimators for $\Gamma^{-1}$ based on $\Gamma_{n}$.

\vspace{.5cm}

The method is theoretically based on the functional calculus for operators (see Dunford and Schwartz, 1988, or Gohberg, Goldberg, Kaashoek, 1991, for instance).

For further purpose we first define the sequence $(\delta_{j})_{j}$ of the smallest differences between distinct eigenvalues of $\Gamma$ as
$$
\left\{\begin{array}[c]{l}
\delta_{1}\ =\ \lambda_{1}-\lambda_{2},\\
\delta_{j}\ =\ \min(\lambda_{j}-\lambda_{j+1},\lambda_{j-1}-\lambda_{j}).
\end{array}\right.
$$
Now take for $(c_{n})_{n\in\mathbb{N}}$ a sequence of strictly positive numbers tending to zero such that $c_n<\lambda_1$ and set
\begin{equation}\label{defkn}
k_{n}=\sup\left\{  p:\lambda_{p}+\delta_{p}/2\geq c_{n}\right\}.
\end{equation}
Then define a class of sequences of positive functions $(f_n)_n$ with support $\left[  c_{n},+\infty\right)$ such that
$$f_{n}\mbox{ is decreasing on }[c_{n},\lambda_{1}+\delta_{1}],\leqno(F.1)$$
\vspace{-1.2cm}
$$\lim_{n\rightarrow+\infty}\sup_{x\geq c_{n}}\left\vert xf_{n}(x)-1\right\vert =0,\leqno(F.2)$$
\vspace{-1.2cm}
$$f_{n}^{\prime}(x)\mbox{ exists for }x\in\lbrack c_{n},+\infty).\leqno(F.3)$$
Moreover, we will make in some cases the additional condition below which will be helpful to reduce the bias of our estimator
$$
\sup_{s\ge c_{n}}\left|  sf_{n}(s)-1\right|  \ =\ o(1/\sqrt{n}).\leqno(H.3)
$$

Now we describe practically the regularization procedure. The eigenvalues of $\Gamma_{n}$ are denoted by $\hat{\lambda}_{j}$ and the associated eigenvectors by $\hat{e}_{j}$. The bounded linear operator $\Gamma_{n}^{\dag}$ is defined the following way

\begin{itemize}
\item Choose a threshold $c_{n}$,
\item Choose a sequence of functions $\left(f_{n}\right)_{n}$ satisfying (F.1)-(F.3),
\item Compute the (functional) PCA of $\Gamma_{n}$ (\textit{i.e.} calculate the eigenvalues $\hat{\lambda}_{j}$ and the eigenvectors $\hat{e}_{j}$),
\item Compute the finite rank operator $\Gamma_{n}^{\dag}$ with the same eigenvectors as $\Gamma_{n}$ and associated eigenvalues $f_{n}(\hat{\lambda}_{j})$ $\left(i.e.\ \Gamma_{n}^{\dag}=\sum_{j=1}^nf_n\left(\hat\lambda_j\right)\hat e_j\otimes\hat e_j\right)$.
\end{itemize}

Obviously $c_{n}$ must be larger than the smallest significatively non-null eigenvalue of $\Gamma_{n}$. Once the threshold $c_n$ and the function $f_n$ (both depending on the sample size $n$) have been chosen, we see that the computation of the estimator of $\rho$ is quite easy through the relation
\begin{equation}
\hat{\rho}\ =\ \Gamma_{n}^{\dag}\Delta_{n}. \label{estimator}
\end{equation}


Now let us give some examples of functions $f_{n}$ and the derived estimators of $\rho$.

\vspace{.5cm}

\noindent\textbf{Example 1.} If $f_{n}(x)=1/x$ when $x\geq c_{n}$ and $0$ elsewhere, condition $(H.3)$ holds and $\Gamma_{n}^{\dag}$ is obtained by simple spectral truncation with threshold $c_{n}.$ The operator $\Gamma_{n}^{\dag}\Gamma_{n}$ is nothing but the projection on a finite dimensional space. Note however that the random dimension of this space, say $d_{n}$, is not necessarily equal to $k_{n}$ (see (\ref{defkn})): for instance we may be in the situation where $\hat{\lambda}_{k_{n}+1}>c_{n}$ and then $d_{n}\geq k_{n}+1$. Unlike $d_n$, $k_{n}$ is non random and was introduced because, as will be seen in the proofs, $P\left(d_{n}\neq k_{n}\right)$ tends to zero fast enough to consider essentially the situation when $d_{n}=k_{n}$. In other words the derived estimator for $\rho$ is asymptotically equivalent to the one considered in Cardot \textit{et al.} (1999).

\vspace{.5cm}

\noindent\textbf{Example 2.} Let $\alpha_{n}$ be some scalar parameter. If $f_{n}(x)=1/(x+\alpha_{n})$ when $x\ge c_{n}$ and 0 elsewhere, we get a ridge-type estimator. Condition (H.3) is satisfied whenever $\alpha_{n}\sqrt{n}/c_{n}\longrightarrow0$.

\vspace{.5cm}

\noindent\textbf{Example 3.} Let $\alpha_{n}$ be some scalar parameter. If $f_{n}(x)=x/(x^{2}+\alpha_{n})$ on its support, $\Gamma^{\dag}_{n}$ is nothing but the Tikhonov regularization of $\Gamma_{n}$. Once more (H.3) holds if $\alpha_{n}\sqrt{n}/c_{n}^{2}\longrightarrow0$.

\vspace{.5cm}

We may define as well, following Mas (1999), a class of approximate for $\Gamma_{n}$ introducing $f_{n,p}(x)=x^{p}/(x+\alpha_{n})^{p+1}$ or $f_{n,p}(x)=x^{p}/(x^{p+1}+\alpha_{n})$, where again $\alpha_{n}$ is some scalar parameter and $p$ some integer.

This procedure is quite general to define regularized version or pseudo inverses for $\Gamma_{n}$. Up to the authors knowledge, all standard techniques for regularizing ill-conditioned matrices or unbounded operators stem from the above functional calculus.


\section{Asymptotic results}


In this section, we mainly announce weak convergence results for the statistical predictor of $Y_{n+1}$ for a new value $X_{n+1}$ obtained by means of estimator defined in (\ref{estimator}), namely $\widehat{Y}_{n+1}=\langle\hat\rho,X_{n+1}\rangle$. Hence, we should study stochastic convergence of
\begin{equation}\label{predictor}
\langle\hat{\rho},X_{n+1}\rangle\,-\,\langle\rho,X_{n+1}\rangle.
\end{equation}
We also look at prediction for a given value of $x\in H$ and study the stochastic convergence of
\begin{equation}
\langle\hat{\rho},x\rangle\,-\,\langle\rho,x\rangle.
\end{equation}
It is important to note that all the results are obtained without assuming any prior knowledge for the rate of decay of the eigenvalues $\lambda_{j}$ of $\Gamma$ to zero. We will see that unfortunately a bias term appears which cannot be removed without very specific assumptions on the sequence on the spectrum of $\Gamma$ and on the smoothness properties of $\rho$.

We begin to investigate the weak convergence for the norm topology on $H$ for our estimate. The next and important result underlines the limits of the functional approach. It tells us that it is not possible to get a general result that would allow to build confidence sets in the functional setting. This highlights the fact that when considering functional data one must take care and multivariate classical results are not necessarily true anymore.

\begin{theorem}
\label{contrex}It is impossible for $\widehat{\rho}-\rho$ to converge in distribution to a non-degenerate r.e. in the norm topology of $H$.
\end{theorem}

The proof of the Theorem is postponed to the end of section \ref{proofs}: it is shown actually that for any normalizing sequence $\alpha_{n}\uparrow+\infty,$ $\alpha_{n}\left(  \widehat{\rho}-\rho\right)$ does not converge in distribution for the norm topology but to a degenerate random element.

Nevertheless this negative result does not mean that it is not possible to get some confidence sets. We have to consider a weak topology (with respect to the inner product), that is to say point-wise confidence bands, and study separately the cases of deterministic and random points. We first give results for the prediction approach.

We define $\Gamma^{\dag}$ as $f_{n}\left(  \Gamma\right)  $. It is important to note that $\Gamma^{\dag}$ depends on the sample size $n$ through the sequence $k_{n}$. From this we take in the following
$$
t_{n,x}\ =\ \sqrt{\sum_{j=1}^{k_{n}}\lambda_{j}\left[  f_{n}(\lambda_{j})\right]  ^{2}\langle x,e_{j}\rangle^{2}}\ =\ \sqrt{\left\Vert\Gamma^{1/2}\Gamma^{\dag}x\right\Vert ^{2}},\ x\in H,
$$
$$
s_{n}\ =\ \sqrt{\sum_{j=1}^{k_{n}}\left[  \lambda_{j}f_{n}(\lambda_{j})\right]  ^{2}}\ =\ \sqrt{tr(\Gamma^{\dag}\Gamma)},
$$
and $\hat{t}_{n,x}$ and $\hat{s}_{n}$ their empirical counterparts based on the $\hat{\lambda}_{j}$'s. Note that the sequence $t_{n,x}$ may either converge or diverge depending on whether $\sum_{j=1}^{+\infty}\lambda_{j}^{-1}\langle x,e_{j}\rangle^{2}=\left\Vert \Gamma^{-1/2}x\right\Vert ^{2}$ is finite or not (\textit{i.e.} whether $x$ is in the range of $\Gamma^{-1/2}$ or not). At the opposite, the term $s_{n}$ always tends to infinity.


\subsection{Weak convergence for the predictor}


We state a weak convergence theorem for the predictor given in (\ref{predictor}). We denote by $\Pi_{k_{n}}$ the projector onto the eigenspace associated to the $k_{n}$ first eigenvalues, and by $\widehat{\Pi}_{k_{n}}$ its empirical counterpart i.e. the projector on the eigenspace associated to $\widehat{\lambda}_{1},\widehat{\lambda}_{2},...,\widehat{\lambda}_{k_{n}}.$\\ 
Assumptions $(H.1)-(H.3)$ are truly basic. They just ensure that the statistical problem is correctly posed. In order to get deep asymptotic results we introduce extra assumptions denoted $(A.1)-(A.3)$.
$$
\sum_{l=1}^{+\infty}\left\vert \left\langle \rho,e_{l}\right\rangle
\right\vert <+\infty.\leqno(A.1)
$$
$$
\mbox{There exists a convex positive function $\lambda,$ such that, at least for $j$ large, $\lambda_{j}=\lambda\left(  j\right)  $.}\leqno(A.2)
$$
We recall the Karhunen-Loève expansion of $X$, that is
$$
X=\sum_{l=1}^{+\infty}\sqrt{\lambda_{l}}\xi_{l}e_{l},
$$
where the $\xi_{l}$'s are centered r.r.v such that $E\xi_{l}\xi_{l^{\prime}}=1$ if $l=l^{\prime}$ and $0$ otherwise. We assume the following assumption for variables $\xi_{l}$
$$
\sup_{l}\mathbb{E}\xi_{l}^{4}\ \leq\ M\ <\ +\infty.\leqno(A.3)
$$

\begin{remark}
Assumption $(A.2)$ is clearly unrestrictive since it holds for standard rates of decrease for the eigenvalues, polynomial or exponential. It implies that
$$
\delta_{k}=\min\left(  \lambda_{k}-\lambda_{k+1},\lambda_{k-1}-\lambda
_{k}\right)  =\lambda_{k}-\lambda_{k+1}.
$$
\end{remark}

\begin{remark}
Simple calculations show that assumption $(A.3)$ implies assumption $(H.1)$, namely that $\mathbb{E}\left\Vert X\right\Vert ^{4}<+\infty$ and does not require any condition on the stochastic dependence within the $\xi_{l}$'s. Besides $(A.3)$ holds for a very large class of real-valued random variables (remind that the $\xi_{l}$'s are subject to $\mathbb{E}\xi_{l}=0$ and $\mathbb{E}\xi_{l}^{2}=1$).
\end{remark}

\begin{theorem}
\label{best}When assumptions $(H.2)-(H.3)$ and $(A.1)-(A.3)$ hold and if
$$
\dfrac{k_{n}^{5/2}\left(  \log k_{n}\right)  ^{2}}{\sqrt{n}}\rightarrow0,
$$
then
$$
\dfrac{\sqrt{n}}{s_{n}}\left(  \left\langle \widehat{\rho},X_{n+1}\right\rangle -\left\langle \Pi_{k_{n}}\rho,X_{n+1}\right\rangle \right)\ \overset{w}{\rightarrow}\ N\left(  0,\sigma_{\varepsilon}^{2}\right).
$$
\end{theorem}

\begin{corollary}
\label{bias}If either $\sup_{p}\left(  \left\vert \left\langle \rho,e_{p}\right\rangle \right\vert p^{5/2}\right)  <+\infty$ or if $\sup_{p}\left(  p^{4}\lambda_{p}\right)  <+\infty$, the bias term $\left\langle\Pi_{k_{n}}\rho,X_{n+1}\right\rangle $ in the previous theorem can be replaced with $\left\langle \rho,X_{n+1}\right\rangle $.
\end{corollary}

\begin{remark}
The term $s_{n}$ always tends to infinity and hence we cannot obtain a "parametric" rate of decay in probability. Besides $s_{n}$ depends on the unknown eigenvalues. It is worth trying to get an "adaptive" version of the above result: replacing the $\lambda_{i}$'s with the $\widehat{\lambda}_{i}$'s leads to a new result with both a random bias and a random normalization term.
\end{remark}

\begin{corollary}
\label{Cor2}The adaptive version of the previous Theorem is
$$
\dfrac{\sqrt{n}}{\widehat{s}_{n}\sigma_{\varepsilon}}\left(  \left\langle
\widehat{\rho},X_{n+1}\right\rangle -\left\langle \Pi_{k_{n}}\rho
,X_{n+1}\right\rangle \right)  \overset{w}{\rightarrow}N\left(  0,1\right)
$$
where
$$
\widehat{s}_{n}=\sqrt{\sum_{j=1}^{k_{n}}\left[  \widehat{\lambda}_{j}f_{n}\left(  \widehat{\lambda}_{j}\right)  \right]  ^{2}}.
$$
\end{corollary}

\begin{remark}
In all the previous results, the variance of the white noise $\sigma_{\varepsilon}^{2}$ is unknown. Replacing $\sigma_{\varepsilon}$ with a convergent estimate of $\sigma_{\varepsilon}$ does not change the Theorems.
\end{remark}


\subsection{Weak convergence for the estimate of $\rho.$}


We are now giving weak convergence results for the prediction at a given value $x$ in $H$.

\begin{theorem}
\label{TH0}Fix any $x$ in $H$. When the assumptions of Theorem \ref{best} hold and if
$$
\sup_{p}\dfrac{\left\vert \left\langle x,e_{p}\right\rangle \right\vert ^{2}}{\lambda_{p}}<+\infty\mathrm{\ and\ }\dfrac{k_{n}^{3}\left(  \log k_{n}\right)  ^{2}}{t_{n,x}\sqrt{n}}\rightarrow0,
$$
then
$$
\dfrac{\sqrt{n}}{t_{n,x}\sigma_{\varepsilon}}\left(  \left\langle\widehat{\rho},x\right\rangle -\left\langle \widehat{\Pi}_{k_{n}}\rho,x\right\rangle \right)  \overset{w}{\rightarrow}N\left(  0,1\right).
$$
\end{theorem}

\begin{remark}
The bias term here is random. It can be seen from the proof of the crucial Proposition \ref{ks} that the situation cannot be improved without very specific (maybe artificial) assumptions either on $\rho$ or on the $\lambda_{i}$'s.
\end{remark}

The normalizing sequence $\dfrac{\sqrt{n}}{t_{n,x}}$ depends on the unknown $\lambda_{j}$'s. It is worth trying to get again an adaptive version of the above theorem (i.e. replace $t_{n,x}$ with $\widehat{t}_{n,x}=\sqrt{\sum_{j=1}^{k_{n}}\widehat{\lambda}_{j}\left[  f_{n}\left(  \widehat{\lambda}_{j}\right)  \right]  ^{2}\left\langle x,\widehat{e}_{j}\right\rangle ^{2}}$).

\begin{corollary}
\label{Cor1}Theorem \ref{TH0} still holds if $t_{n,x}$ is replaced with its empirical counterpart $\widehat{t}_{n,x}.$
\end{corollary}

The following Remark is crucial since it brings out once more what seems to be a typical feature of the functional setting.

\begin{remark}
As seen before the sequence $t_{n,x}$ may either converge or diverge. Indeed, if $\left\Vert \Gamma^{-1/2}x\right\Vert $ is finite the normalization sequence grows surprisingly at a parametric rate (i.e $\sqrt{n}$). This could be understood as an extra-smoothing of the estimate $\widehat{\rho}$ through the integrals involving the scalar product. But in terms of prediction this fact could be misleading. This "extra-smoothing" is indeed an excessive and artificial smoothing since $P\left(  \left\Vert \Gamma^{-1/2}X_{n+1}\right\Vert<+\infty\right)  =0.$ This also means the realizations of $X$ do not belong with probability one to the reproducing kernel Hilbert space associated to its covariance function (Hajek, 1962). In other words the results of this section are given for the sake of completeness and to explore the analytical properties of our estimates. For these reasons and if prediction is under concern, only $\left\langle \widehat{\rho},X_{n+1}\right\rangle $ should be considered and studied. In a multivariate setting all these considerations make no sense, since the situation is simpler (in fact, usually $P\left(  \left\Vert\Gamma^{-1/2}X_{n+1}\right\Vert <+\infty\right)  =1$ because $\Gamma^{-1/2}$ is bounded when $\Gamma$ is a full rank covariance matrix).
\end{remark}

Within the proofs it is readily seen that assumption $\left(  A.1\right)$ plays a crucial role in getting a non random bias term. The next Proposition illustrates this situation.

\begin{proposition}
\label{contrex2}Assume that $\lambda_{k}=k^{-1-\alpha},$ that $\left\langle x,e_{k}\right\rangle ^{2}=k^{-1-\beta}$ with $\beta>1+\alpha.$ Then if $\sum_{j=1}^{+\infty}j^{1-\beta}\left\langle \rho,e_{j}\right\rangle^{2}=+\infty$, the sequence $\dfrac{\sqrt{n}}{t_{n,x}\sigma_{\varepsilon}}\left\langle \left(  \widehat{\Pi}_{k_{n}}-\Pi_{k_{n}}\right)  \rho,x\right\rangle $ may not be bounded in probability even if the random variables $X_{i}$, $i=1,\ldots,n$, are i.i.d. centered, Gaussian.
\end{proposition}

\begin{remark}
The condition $\beta>1+\alpha$ just ensures that $\left\Vert \Gamma^{-1/2}x\right\Vert <+\infty.$ Besides if $\left\langle \rho,e_{j}\right\rangle ^{2}=j^{-1-\gamma},$ $\sum_{j=1}^{n}j^{1-\beta}\left\langle\rho,e_{j}\right\rangle ^{2}$ diverges whenever $\beta+\gamma<1$ which implies that $\sum_{j=1}^{+\infty}\left\vert \left\langle \rho,e_{j}\right\rangle \right\vert =+\infty.$
\end{remark}

In fact the assumption on the location of $\rho$ mentioned in the Proposition should be understood as smoothness conditions.


\section{Concluding remarks}


One important application of previous results is the construction of confidence sets for prediction.
In real life problems, the regression function $\rho$ is unknown but Corollary \ref{bias} allows us to build confidence sets.
Let $q_{\alpha}$ be the quantile of order $1-\alpha/2 $ of a Gaussian random variable with mean $0$ and unit variance, we get under previous assumptions the following confidence set for prediction,
\begin{eqnarray}
\lim_{n \to \infty} P \left(  \frac{\sqrt{n}}{\hat{\sigma} \ \widehat{s}_{n}} \left\vert<\widehat{\rho},X_{n+1}> - <\rho,X_{n+1}> \right\vert \geq q_{\alpha}\right) \ = \ 1 - \alpha \ .
\end{eqnarray}
A simulation study (Cardot \textit{et al.}, 2004) has shown that such confidence sets are accurate even for moderate sample sizes, \textit{i.e.} for $n$ around 100.

From a mathematical points of view, one of the main novelty of this work relies on the facts that no prior information on the eigenvalues is assumed and the dimension sequence $k_{n}$ does not depend on the rate of decrease of these eigenvalues. As a consequence $k_{n}$ increase rather slowly, but not that much for a non parametric model. 
 Nevetheless, let us  notice that  this situation may be significantly improved if some information on the eigenvalues is available. 

From Theorem \ref{best} it is possible to derive a general bound for the $L^{2}$ prediction error. Simple calculations (see the proof of Corollary \ref{bias}) lead to :
\begin{equation}
\left\langle \widehat{\rho}-\rho,X_{n+1}\right\rangle ^{2}=O_{\mathbb{P}}\left(  \frac{s_{n}}{n}\right)  +O_{\mathbb{P}}\left(  \sum_{j=k_{n}+1}^{\infty}\lambda_{j}\langle\rho,e_{j}\rangle^{2}\right).\label{biasboundKL}
\end{equation}
Thus, it is not possible to go further without imposing more precise hypotheses on the smoothness of function $\rho$ with respect to the basis of eigenfunctions $e_{j}$ and the rate of decay of the eigenvalues $\lambda_{j}$ as remarked sooner in the article. Nevertheless, it was seen that the second term on the right in (\ref{biasboundKL}) can converge rapidly to zero in some situations. Besides assumption $(A.1)$ provides us with some kind of uniformity with respect to $\rho$ when the latter belongs to a subset of $H.$ Naturally, with these remarks we have in mind the study of the minimax  rate of $L^{2}$ risk for the class of our predictor.


\section{Proofs\label{proofs}}


Along the proofs we suppose that $(H.1)-(H.3)$ hold. The letter $C$ will always stand for any (nonrandom and universal) positive constant. For any bounded operator $T$ defined and with values in $H$ we classically set
$$
\left\Vert T\right\Vert _{\infty}=\sup_{x\in B_{1}}\left\Vert Tx\right\Vert,
$$
where $B_{1}$ is the unit ball of $H.$ We will quite often make use of the following facts.

\begin{itemize}
\item For any $u$ in $H,$
$$
\mathbb{E}\left\langle X,u\right\rangle ^{2}=\left\langle \Gamma
u,u\right\rangle =\left\Vert \Gamma^{1/2}u\right\Vert ^{2}.
$$
\item For a sufficiently large $i$, $\lambda_{i}\leq\dfrac{C}{i\log i}.$
\item The Hilbert-Schmidt norm is more precise than the classical norm for operators. Hence if $T$ is Hilbert-Schmidt
$$
\left\Vert T\right\Vert _{\infty}\leq\left\Vert T\right\Vert _{HS}=\sqrt{\sum_{p}\left\Vert Tu_{p}\right\Vert ^{2}},
$$
where $\left(  u_{p}\right)  _{p\in\mathbb{N}}$ is any complete orthonormal sequence in $H$.
\end{itemize}
From definitions of $\hat{\rho}$ and $U_{n}$ we have
$$
\widehat{\rho}=\Gamma_{n}^{\dagger}\Gamma_{n}\rho+\left(  \Gamma_{n}^{\dagger}-\Gamma^{\dagger}\right)  U_{n}+\Gamma^{\dagger}U_{n},
$$
from which the forthcoming decomposition is trivial

\label{decom}
\begin{equation}
\widehat{\rho}-\Pi_{k_{n}}\rho=T_{n}+S_{n}+R_{n}+Y_{n}, \label{decomp}
\end{equation}
where
$$
T_{n}  =\left(  \Gamma_{n}^{\dagger}\Gamma_{n}-\widehat{\Pi}_{k_{n}}\right)  \rho,\qquad S_{n}=\left(  \Gamma_{n}^{\dagger}-\Gamma^{\dagger}\right)  U_{n},\qquad R_{n}   =\Gamma^{\dagger}U_{n}\qquad Y_{n}=\left(  \widehat{\Pi}_{k_{n}}-\Pi_{k_{n}}\right)  \rho.
$$
We also denote
$$
L_{n}=\Pi_{k_{n}}\rho-\rho.
$$
The proofs are tiled into four subsections. After a brief introduction on operator-valued analytic functions, we begin with providing useful convexity inequalities for the eigenvalues and subsequent moment bounds. The second part shows that all the bias terms but $L_{n},$ say $T_{n},$ $S_{n}$ and $Y_{n}$ tend to zero in probability when correctly normalized. Weak convergence of $R_{n}$ is proved in the short third subsection. The last part provides the main results of the paper by collecting the Lemmas and Propositions previously proved.


\subsection{Preliminary results}


All along the proofs we will need auxiliary results from perturbation theory for bounded operators. It is of much help to have basic notions about spectral representation of bounded operators and perturbation theory. We refer to Dunford-Schwartz (1988, Chapter VII.3) or to Gohberg, Goldberg and Kaashoek (1991) for an introduction to functional calculus for operators related with Riesz integrals. 

Let us denote by $\mathcal{B}_{i}$ the oriented circle of the complex plane with center $\lambda_{i}$ and radius $\delta_{i}/2$ and define
$$
\mathcal{C}_{n}=\bigcup_{i=1}^{k_{n}}\mathcal{B}_{i}\ .
$$
The open domain whose boundary is $\mathcal{C}_{n}$ is not connected but however we can apply the functional calculus for bounded operators (see Dunford-Schwartz Section VII.3 Definitions 8 and 9). We also need to change slightly the definition of the sequence of functions $\left(f_{n}\right)_n$ by extending it to the complex plane, more precisely to $\mathcal{C}_{n}.$ We admit that it is possible to extend $f_{n}$ to an analytic function $\widetilde{f_{n}}$ defined on the interior of $\mathcal{C}_{n}$ (in the plane) such that $\sup_{z\in\mathcal{C}_{n}}\left\vert \widetilde{f}_{n}\left(  z\right)\right\vert \leq C\sup_{x\in\left[  c_{n},\lambda_{1}+\delta_1\right]  }\left\vert f_{n}\left(  x\right)  \right\vert$. For instance if $f_{n}\left(  x\right)=(1/x)1\!\!1_{[c_n,+\infty)}(x)$, take $\widetilde{f}_{n}\left(  z\right)  =(1/z)1\!\!1_{\mathcal{C}_{n}}(z)$. Results from perturbation theory yield
\begin{align}
\Pi_{k_{n}}  &  =\dfrac{1}{2\pi\iota}\int_{\mathcal{C}_{n}}\left(z-\Gamma\right)  ^{-1}dz,\label{respro}\\
\Gamma^{\dag}  &  =\dfrac{1}{2\pi\iota}\int_{\mathcal{C}_{n}}\left(z-\Gamma\right)  ^{-1}f_{n}\left(  z\right)  dz, \label{resinv}
\end{align}
where $\iota^{2}=-1.$

We introduce also the square root of symmetric operators: if $T$ is a positive self-adjoint operator (random or not), we denote by $\left(  zI-T\right)^{1/2}$ the symmetric operator whose eigenvectors are the same as $T$ and whose eigenvalues are the complex square root of $z-\lambda_{k}$, $k\in\mathbb{N}$, denoted $\left(  z-\lambda_{k}\right)  ^{1/2}.$

\begin{lemma}
\label{trick1}Consider two positive integers $j$ and $k$ large enough and such that $k>j$. Then
\begin{equation}
j\lambda_{j}\ \geq\ k\lambda_{k}\qquad\mbox{and}\qquad\lambda_{j}-\lambda_{k}\geq\left(  1-\dfrac{j}{k}\right)  \lambda_{j}.
\label{t1}
\end{equation}
Besides
\begin{equation}
\sum_{j\geq k}\lambda_{j}\leq\left(  k+1\right)  \lambda_{k} \label{t2}.
\end{equation}
\end{lemma}

\begin{proof}
We set for notational convenience $\lambda_{j}=\varphi\left(  1/j\right)$ where $\varphi$ is, by assumption $(A.2)$, a convex function defined on the interval $[0,1]$ such that that $\varphi\left(0\right)  =0$ and $\varphi(1)=\lambda_{1}$.\\
The two inequalities in (\ref{t1}) follows directly from the well known inequalities for convex functions
$$
\frac{\varphi(x_{1})-\varphi(x_{0})}{x_{1}-x_{0}}\leq\frac{\varphi
(x_{2})-\varphi(x_{0})}{x_{2}-x_{0}}\leq\frac{\varphi(x_{2})-\varphi(x_{1})}{x_{2}-x_{1}}\ ,\ 0\leq x_{0}<x_{1}<x_{2}\leq1,
$$
and by taking $x_{0}=0,\ x_{1}=1/k$ and $x_{2}=1/j$.\\ 
Set $\mu_{k}=\sum_{l\geq k}\lambda_{l}.$ It is easy to see that the sequence $\left(\mu_{k}\right)_k$ satisfies assumption $(A.2)$. Indeed for all $k$
$$
\mu_{k}-\mu_{k+1}\leq\mu_{k-1}-\mu_{k},
$$
which is a sufficient condition to construct a convex function $\mu\left(k\right)=\mu_k$. We can then apply the second part of (\ref{t1}) with $\mu_{k+1}$ instead of $\lambda_k$ and $\mu_k$ instead of $\lambda_j$, which yields
$$
\mu_{k}-\mu_{k+1}=\lambda_{k}\geq\dfrac{1}{k+1}\mu_{k},
$$
and (\ref{t2}) is proved.
\end{proof}

\begin{lemma}
\label{trick2}The following is true for $j$ large enough
$$
\sum_{l\neq j}\dfrac{\lambda_{l}}{\left\vert \lambda_{l}-\lambda_{j}\right\vert }\leq Cj\log j.
$$

\end{lemma}

\begin{proof}
We are first going to decompose the sum into three terms
$$
\sum_{l\neq j}\dfrac{\lambda_{l}}{\left\vert \lambda_{l}-\lambda
_{j}\right\vert }=\mathcal{T}_{1}+\mathcal{T}_{2}+\mathcal{T}_{3},
$$
where
$$
\mathcal{T}_{1}=\sum_{l=1}^{j-1}\dfrac{\lambda_{l}}{\lambda_{l}-\lambda_{j}},\ \mathcal{T}_{2}=\sum_{l=j+1}^{2j}\dfrac{\lambda_{l}}{\lambda_{j}-\lambda_{l}},\ \mathcal{T}_{3}=\sum_{l=2j+1}^{+\infty}\dfrac{\lambda_{l}}{\lambda_{j}-\lambda_{l}}.
$$
Applying Lemma \ref{trick1} we get
$$
\mathcal{T}_{1}=\sum_{l=1}^{j-1}\dfrac{\lambda_{l}}{\lambda_{l}-\lambda_{j}}\leq j\sum_{l=1}^{j-1}\dfrac{1}{j-l}\leq C_{1}j\log j,
$$
where $C_{1}$ is some positive constant. Also, applying once more (\ref{t1}) then (\ref{t2}), we get
\begin{align*}
\mathcal{T}_{2}  &  =\sum_{l=j+1}^{2j}\dfrac{\lambda_{l}}{\lambda_{j}-\lambda_{l}}\leq\sum_{l=j+1}^{2j}\dfrac{\lambda_{l}}{\lambda_{j}}\dfrac{l}{l-j}\\
&  \leq2j\sum_{l=j+1}^{2j}\dfrac{1}{l-j}\leq C_{2}j\log j,
\end{align*}
and
$$
\mathcal{T}_{3}\leq\sum_{l=2j+1}^{+\infty}\dfrac{\lambda_{l}}{\lambda_{j}-\lambda_{l}}\leq\dfrac{\sum_{l=2j+1}^{+\infty}\lambda_{l}}{\lambda_{j}-\lambda_{2j}}\leq2\dfrac{\sum_{l=2j+1}^{+\infty}\lambda_{l}}{\lambda_{j}}\leq C_{3}j.
$$
Hence the result follows and Lemma \ref{trick2} is proved.
\end{proof}

\begin{lemma}
\label{thal} We have for $j$ large enough
\begin{equation}
\mathbb{E}\sup_{z\in\mathcal{B}_{j}}\left\Vert \left(  zI-\Gamma\right)^{-1/2}\left(  \Gamma_{n}-\Gamma\right)  \left(  zI-\Gamma\right)^{-1/2}\right\Vert _{\infty}^{2}  \leq\dfrac{C}{n}\left(  j\log j\right)^{2},\label{majcov}
\end{equation}
and
$$
\mathbb{E}\sup_{z\in\mathcal{B}_{j}}\left\Vert \left(  zI-\Gamma\right)
^{-1/2}X_{1}\right\Vert ^{2}   \leq Cj\log j.
$$
\end{lemma}

\begin{proof}
Take $z\in\mathcal{B}_{j}.$ By bounding the sup norm by the Hilbert-Schmidt one (see above), we get
\begin{align*}
&  \left\Vert \left(  zI-\Gamma\right)  ^{-1/2}\left(  \Gamma_{n}-\Gamma\right)  \left(  zI-\Gamma\right)  ^{-1/2}\right\Vert _{\infty}^{2}\\
&  \leq \sum_{l=1}^{+\infty}\sum_{k=1}^{+\infty}\left\langle \left(zI-\Gamma\right)  ^{-1/2}\left(  \Gamma_{n}-\Gamma\right)  \left(zI-\Gamma\right)  ^{-1/2}\left(  e_{l}\right)  ,e_{k}\right\rangle ^{2}\\
&  \leq\sum_{l,k=1}^{+\infty}\dfrac{\left\langle \left(  \Gamma_{n}-\Gamma\right)  \left(  e_{l}\right)  ,e_{k}\right\rangle ^{2}}{\left\vert z-\lambda_{l}\right\vert \left\vert z-\lambda_{k}\right\vert }\\
&  \leq4\sum_{\substack{l,k=1,\\l,k\neq j}}^{+\infty}\dfrac{\left\langle\left(  \Gamma_{n}-\Gamma\right)  \left(  e_{l}\right)  ,e_{k}\right\rangle^{2}}{\left\vert \lambda_{j}-\lambda_{l}\right\vert \left\vert \lambda_{j}-\lambda_{k}\right\vert }+2\sum_{\substack{k=1,\\k\neq j}}^{+\infty}\dfrac{\left\langle \left(  \Gamma_{n}-\Gamma\right)  \left(  e_{j}\right),e_{k}\right\rangle ^{2}}{\delta_{j}\left\vert z-\lambda_{k}\right\vert}+\dfrac{\left\langle \left(  \Gamma_{n}-\Gamma\right)  \left(  e_{j}\right),e_{j}\right\rangle ^{2}}{\delta_{j}^{2}},
\end{align*}
since it can be checked that whenever $z=\lambda_{j}+\dfrac{\delta_{j}}{2}e^{\iota\theta}\in\mathcal{B}_{j}$ and $i\neq j$
$$
\left\vert z-\lambda_{i}\right\vert =\left\vert \lambda_{j}-\lambda_{i}+\dfrac{\delta_{j}}{2}e^{\iota\theta}\right\vert \geq\left\vert \lambda_{j}-\lambda_{i}\right\vert -\dfrac{\delta_{j}}{2}\geq\left\vert \lambda_{j}-\lambda_{i}\right\vert /2.
$$
Besides
$$
\mathbb{E}\left\langle \left(  \Gamma_{n}-\Gamma\right)  \left(e_{l}\right)  ,e_{k}\right\rangle ^{2}\ =\ \dfrac{1}{n}\left[  \mathbb{E}\left(  \left\langle X_{1},e_{k}\right\rangle ^{2}\left\langle X_{1},e_{l}\right\rangle ^{2}\right)-\left\langle \Gamma\left(  e_{l}\right)  ,e_{k}\right\rangle ^{2}\right]\leq\dfrac{M}{n}\lambda_{l}\lambda_{k},
$$
when assumption $(A.3)$ holds. Finally
\begin{align*}
&  \mathbb{E}\sup_{z\in\mathcal{B}_{j}}\left\Vert \left(  zI-\Gamma\right)^{-1/2}\left(  \Gamma_{n}-\Gamma\right)  \left(  zI-\Gamma\right)
^{-1/2}\right\Vert _{\infty}^{2}\\
&  \leq\dfrac{M}{n}\left[  \sum_{\substack{l,k=1,\\l,k\neq j}}^{+\infty}\dfrac{\lambda_{l}\lambda_{k}}{\left\vert \lambda_{j}-\lambda_{l}\right\vert\left\vert \lambda_{j}-\lambda_{k}\right\vert }+\dfrac{\lambda_{j}}{\delta_{j}}\sum_{k=1,k\neq j}^{+\infty}\dfrac{\lambda_{k}}{\left\vert \lambda_{j}-\lambda_{l}\right\vert }+\left(  \dfrac{\lambda_{j}}{\delta_{j}}\right)^{2}\right] \\
&  =\dfrac{M}{n}\left[  \left(  \sum_{k=1,k\neq k}^{+\infty}\dfrac{\lambda_{k}}{\left\vert \lambda_{j}-\lambda_{k}\right\vert }\right)  ^{2}%
+\dfrac{\lambda_{j}}{\delta_{j}}\sum_{k=1,k\neq j}^{+\infty}\dfrac{\lambda_{k}}{\left\vert \lambda_{j}-\lambda_{l}\right\vert }+\left(  \dfrac
{\lambda_{j}}{\delta_{j}}\right)  ^{2}\right]  .
\end{align*}
It suffices now to apply Lemmas \ref{trick1} and \ref{trick2} to get the desired result. The same method leads to proving the second part of the display. Lemma \ref{thal} is proved.
\end{proof}

\begin{lemma}
\label{you}Denoting
$$
\mathcal{E}_{j}\left(  z\right)  =\left\{  \left\Vert \left(  zI-\Gamma\right)  ^{-1/2}\left(  \Gamma_{n}-\Gamma\right)  \left(  zI-\Gamma\right)^{-1/2}\right\Vert _{\infty}<1/2,z\in\mathcal{B}_{j}\right\},
$$
The following holds
$$
\left\Vert \left(  zI-\Gamma\right)  ^{1/2}\left(  zI-\Gamma_{n}\right)^{-1}\left(  zI-\Gamma\right)  ^{1/2}\right\Vert _{\infty}1\!\!1_{\mathcal{E}_{j}\left(  z\right)  }\leq C,\quad a.s.
$$
where $C$ is some positive constant. Besides
\begin{equation}
\mathbb{P}\left(  \mathcal{E}_{j}^{c}\left(  z\right)  \right)  \leq\dfrac{j\log j}{\sqrt{n}}. \label{anytime}
\end{equation}
\end{lemma}

\begin{proof}
We have successively
$$
\left(  zI-\Gamma_{n}\right)  ^{-1}=\left(  zI-\Gamma\right)  ^{-1}+\left(zI-\Gamma\right)  ^{-1}\left(  \Gamma-\Gamma_{n}\right)  \left(  zI-\Gamma
_{n}\right)  ^{-1},
$$
hence
\begin{equation}
\left(  zI-\Gamma\right)  ^{1/2}\left(  zI-\Gamma_{n}\right)  ^{-1}\left(zI-\Gamma\right)  ^{1/2}=I+\left(  zI-\Gamma\right)  ^{-1/2}\left(  \Gamma-\Gamma_{n}\right)\left(  zI-\Gamma_{n}\right)  ^{-1}\left(  zI-\Gamma\right)  ^{1/2},
\label{invborn}
\end{equation}
and
\begin{equation}
\left[  I+\left(  zI-\Gamma\right)  ^{-1/2}\left(  \Gamma_{n}-\Gamma\right)\left(  zI-\Gamma\right)  ^{-1/2}\right]  \left(  zI-\Gamma\right)^{1/2}\left(  zI-\Gamma_{n}\right)  ^{-1}\left(  zI-\Gamma\right)  ^{1/2}=I.
\label{invborn2}
\end{equation}
It is a well known fact that if the linear operator $T$ satisfies $\left\Vert T\right\Vert _{\infty}<1$ then $I+T$ is an invertible and its inverse is bounded and given by formula
$$
\left(  I+T\right)  ^{-1}=I-T+T^{2}-...
$$
From (\ref{invborn}) and (\ref{invborn2}) we deduce that
\begin{align*}
&  \left\|  \left(  zI-\Gamma\right)  ^{1/2}\left(  zI-\Gamma_{n}\right)^{-1}\left(  zI-\Gamma\right)  ^{1/2}\right\|_\infty  1\!\!1_{\mathcal{E}_{j}\left(z\right)  }\\
&  =\left\|\left[  I+\left(  zI-\Gamma\right)  ^{-1/2}\left(  \Gamma_{n}-\Gamma\right)  \left(  zI-\Gamma\right)  ^{-1/2}\right]  ^{-1}\right\|_\infty1\!\!1_{\mathcal{E}_{j}\left(  z\right)  }\leq2,\quad a.s.
\end{align*}
Now, the bound in (\ref{anytime}) stems easily from Markov inequality and (\ref{majcov}) in Lemma \ref{thal}. This finishes the proof of the Lemma.
\end{proof}

\medskip

The empirical counterparts of (\ref{respro}) and (\ref{resinv}) -mentioned above- involve a random contour, say $\widehat{\mathcal{B}}_{i},$ centered at $\widehat{\lambda}_{i}.$ It should be noted that these contours cannot be replaced by the $\mathcal{B}_{i}$'s since the latter may contain more than $k_{n}$ eigenvalues of $\Gamma_{n}$. The aim of the following Lemma is to find sufficient conditions under which $\widehat{\mathcal{B}}_{i}$ may be replaced with $\mathcal{B}_{i}.$ In other words, we have to check that for a sufficiently large $n$ the $p^{th}$ eigenvalue of $\Gamma_{n}$ is close enough from the $p^{th}$ eigenvalue of $\Gamma.$ Before stating this first lemma, we introduce the following event
$$
\mathcal{A}_{n}=\left\{  \forall j\in\left\{  1,...,k_{n}\right\}|\dfrac{\left\vert \widehat{\lambda}_{j}-\lambda_{j}\right\vert }{\delta_{j}}<1/2\right\}.
$$

\begin{lemma}
\label{restric}If $\dfrac{k_{n}^{2}\log k_{n}}{\sqrt{n}}\rightarrow0$, then
$$
\dfrac{1}{2\pi\iota}\int_{\mathcal{C}_{n}}\left(  z-\Gamma_{n}\right)^{-1}dz=\widehat{\Pi}_{k_{n}}1\!\!1_{\mathcal{A}_{n}}+r_{n},
$$
where $r_{n}$ is a random operator satisfying $\sqrt{n}r_{n}\overset{\mathbb{P}}{\rightarrow}0$ in the operator norm.
\end{lemma}

\begin{proof}
When the event $\mathcal{A}_{n}$ holds, the $k_{n}$ first empirical eigenvalues $\widehat{\lambda}_{j}$ lie in $\mathcal{B}_{j}$ and then
$$
\widehat{\Pi}_{k_{n}}=\dfrac{1}{2\pi\iota}\int_{\widehat{\mathcal{C}}_{n}}\left(  z-\Gamma_{n}\right)  ^{-1}dz=\dfrac{1}{2\pi\iota}\int_{\mathcal{C}_{n}}\left(  z-\Gamma_{n}\right)  ^{-1}dz.
$$
From this it is clear that
$$
\dfrac{1}{2\pi\iota}\int_{\mathcal{C}_{n}}\left(  z-\Gamma_{n}\right)^{-1}dz=\widehat{\Pi}_{k_{n}}1\!\!1_{\mathcal{A}_{n}}+1\!\!1_{\mathcal{A}_{n}^{c}}\dfrac{1}{2\pi\iota}\int_{\mathcal{C}_{n}}\left(  z-\Gamma_{n}\right)  ^{-1}dz.
$$
Denoting $r_{n}=1\!\!1_{\mathcal{A}_{n}^{c}}\dfrac{1}{2\pi\iota}\int_{\mathcal{C}_{n}}\left(  z-\Gamma_{n}\right)  ^{-1}dz$, we see that, since $\left\Vert \dfrac{1}{2\pi\iota}\int_{\mathcal{C}_{n}}\left(z-\Gamma_{n}\right)  ^{-1}dz\right\Vert _{\infty}=1$, we have for $\varepsilon>0$
$$
\mathbb{P}\left(  \sqrt{n}\left\Vert r_{n}\right\Vert _{\infty}>\varepsilon\right)  \leq\mathbb{P}\left(  1\!\!1_{\mathcal{A}_{n}^{c}}>\varepsilon\right)=\mathbb{P}\left(\mathcal{A}_{n}^{c}\right).
$$
It remains to find a bound for $\mathbb{P}\left(\mathcal{A}_{n}^{c}\right)$. We have
\begin{align}   
\mathbb{P}\left(  \mathcal{A}_{n}^{c}\right) & \leq\sum_{j=1}^{k_{n}}\mathbb{P}\left(  \left\vert \widehat{\lambda}_{j}-\lambda_{j}\right\vert>\delta_{j}/2\right) \nonumber\\
&  \leq2\sum_{j=1}^{k_{n}}\dfrac{\mathbb{E}\left\vert \widehat{\lambda}_{j}-\lambda_{j}\right\vert }{\delta_{j}}=\dfrac{2}{\sqrt{n}}\sum_{j=1}^{k_{n}}\dfrac{\lambda_{j}}{\delta_{j}}\dfrac{\sqrt{n}\mathbb{E}\left\vert\widehat{\lambda}_{j}-\lambda_{j}\right\vert }{\lambda_{j}}. 
\label{aa}
\end{align}
In order to get a uniform bound with respect to $j$ of the latter expectation we follow the same arguments as Bosq (2000), proof of Theorem 4.10 p.122-123. In Bosq, the setting is quite more general but however his Theorem 4.10 ensures that in our framework the asymptotic behaviour of $\sqrt{n}\dfrac{\left\vert \widehat{\lambda}_{j}-\lambda_{j}\right\vert}{\lambda_{j}}$ is the same as $\sqrt{n}\dfrac{\left\vert \left\langle \left(\Gamma_{n}-\Gamma\right)  e_{j},e_{j}\right\rangle \right\vert }{\lambda_{j}}$. From assumption (A.3), we get
\begin{equation}
\sqrt{n}\dfrac{\mathbb{E}\left\vert \left\langle \left(  \Gamma_{n}-\Gamma\right)  e_{j},e_{j}\right\rangle \right\vert }{\lambda_{j}}  \leq\dfrac{\sqrt{\mathbb{E}\left\vert \left\langle X_{1}e_{j}\right\rangle^{4}-\lambda_{j}^{2}\right\vert }}{\lambda_{j}}\leq C ,\label{aaa}
\end{equation}
where $C$ does not depend on $j$. 
From (\ref{aa}) and (\ref{aaa}) we deduce, applying Lemma \ref{trick1} once more, that
$$
\mathbb{P}\left(  \mathcal{A}_{n}^{c}\right)  \leq\dfrac{C}{\sqrt{n}}\sum_{j=1}^{k_{n}}\dfrac{\lambda_{j}}{\delta_{j}}\leq\dfrac{C}{\sqrt{n}}\sum_{j=1}^{k_{n}}j\log j\leq\dfrac{C}{\sqrt{n}}k_{n}^{2}\log k_{n},
$$
from which the result follows.
\end{proof}

\noindent It may be easily proved that the same result as in the preceding Lemma holds with $\Gamma_{n}^{\dagger}$ instead of $\widehat{\Pi}_{k_{n}}$. From now on we will implicitly work on the space $\mathcal{A}_{n}$ and then write 
$$\widehat{\Pi}_{k_{n}}\ =\ \left(  \dfrac{1}{2\pi\iota}\int_{\mathcal{C}_{n}}\left(  z-\Gamma_{n}\right)  ^{-1}dz\right),$$
and
$$
\Gamma_{n}^{\dagger}    \ =\ \left(  \dfrac{1}{2\pi\iota}\int_{\mathcal{C}_{n}}\left(  z-\Gamma_{n}\right)  ^{-1}\widetilde{f}_{n}\left(  z\right)dz\right).
$$
We will also abusively denote  $\Pi_{k_{n}}1\!\!1_{\mathcal{A}_{n}}$ by $\Pi_{k_{n}}$ and $\Gamma^{\dagger}1\!\!1_{\mathcal{A}_{n}}$ by  $\Gamma^{\dagger}$.

\begin{remark}
In fact thanks to Lemma \ref{restric}, we can deal with all our random
elements as if almost surely all the random eigenvalues were in their
associated circles $\mathcal{B}_{j}.$ The reader should keep this fact in mind
all along the forthcoming proofs. The condition on $k_{n}$ needed on the Lemma
is clearly weaker that the ones which appear for the main results to hold.
\end{remark}


\subsection{Bias terms}


As announced above this subsection is devoted to the bias terms $S_{n},T_{n}$ and $Y_{n}.$ A bound is also given for $L_{n}$ for further purpose. We first begin with the term $T_{n}$ for which we have the following lemma.

\begin{lemma}
\label{Tn}If $(H.3)$ holds
$$
\left\Vert T_{n}\right\Vert _{\infty}=\left\Vert \left(  \Gamma_{n}^{\dagger}\Gamma_{n}-\widehat{\Pi}_{k_{n}}\right)  \rho\right\Vert _{\infty}=o_{P}\left(  \dfrac{1}{\sqrt{n}}\right).
$$

\end{lemma}

\begin{proof}
Obviously $\Gamma_{n}^{\dagger}\Gamma_{n}-\widehat{\Pi}_{k_{n}}$ is a self-adjoint random operator whose eigenvalues are the $\left(  \widehat{\lambda}_{j}f_{n}\left(  \widehat{\lambda}_{j}\right)  -1\right)  _{1\leq j\leq k_{n}}$ and $0$ otherwise. So we have
$$
\left\Vert \Gamma_{n}^{\dagger}\Gamma_{n}-\widehat{\Pi}_{k_{n}}\rho\right\Vert_{\infty}\leq C\sup_{s\geq c_{n}}\left(  \left\vert sf_{n}\left(  s\right)-1\right\vert \right).
$$

If assumption $(H.3)$ holds, the last term above is an $o\left(  1/\sqrt{n}\right)$, which proves the second equality.
\end{proof}

\begin{lemma}
\label{Ln}The two following bounds are valid
\begin{align*}
\sqrt{\dfrac{n}{k_{n}}}\mathbb{E}\left\vert L_{n}\right\vert  &  \leq \sqrt{\dfrac{n}{k_{n}}}\left\vert \left\langle \rho,e_{k_{n}}\right\rangle\right\vert \sqrt{\sum_{l\ge k_{n}+1}\lambda_{l}},\\
\sqrt{\dfrac{n}{k_{n}}}\mathbb{E}\left\vert L_{n}\right\vert  &  \leq\dfrac{\lambda_{k_{n}}}{k_{n}}\sqrt{\dfrac{n}{\log k_{n}}}\sqrt{\sum_{l\ge k_{n}+1}\left\langle \rho,e_{l}\right\rangle }.
\end{align*}
\end{lemma}

\begin{proof}We have
\begin{align*}
\mathbb{E}\left\vert \left\langle \left(  I-\Pi_{k_{n}}\right)  \rho,X_{n+1}\right\rangle \right\vert  &  \leq\sqrt{\mathbb{E}\sum_{l=k_{n}+1}\left\langle \rho,e_{l}\right\rangle^{2}\left\langle X_{n+1},e_{l}\right\rangle ^{2}}\\
&  =\sqrt{\sum_{l\ge k_{n}+1}\lambda_{l}\left\langle \rho,e_{l}\right\rangle^{2}}\\
&  \leq\left\{
\begin{tabular}[c]{l}
$\left\vert \left\langle \rho,e_{k_{n}}\right\rangle \right\vert \sqrt{\sum_{l\ge k_{n}+1}\lambda_{l}}$\\
$\dfrac{\lambda_{k_{n}}}{\sqrt{k_{n}\log k_{n}}}\sqrt{\sum_{l\ge k_{n}+1}\left\langle \rho,e_{l}\right\rangle }$,
\end{tabular}\right.
\end{align*}
since $\lambda_{l}$ and $\left\vert \left\langle \rho,e_{l}\right\rangle
\right\vert $ are absolutely summing sequences.
\end{proof}

\begin{proposition}
\label{ks}If $\dfrac{1}{\sqrt{n}}k_{n}^{5/2}\left(  \log k_{n}\right)^{2}\rightarrow0$ as $n$ goes to infinity, then
$$
\sqrt{\dfrac{n}{k_{n}}}\left\langle \left(  \widehat{\Pi}_{k_{n}}-\Pi_{k_{n}}\right)  \rho,X_{n+1}\right\rangle \overset{\mathbb{P}}{\rightarrow}0.
$$

\end{proposition}

\begin{proof}
The proof of the Proposition is the keystone of the paper. We begin with
\begin{align*}
\left(  \widehat{\Pi}_{k_{n}}-\Pi_{k_{n}}\right)   &  =\dfrac{1}{2\pi\iota}\sum_{j=1}^{k_{n}}\int_{\mathcal{B}_{j}}\left[  \left(  zI-\Gamma_{n}\right)^{-1}-\left(  zI-\Gamma\right)  ^{-1}\right]  dz\\
&  =\dfrac{1}{2\pi\iota}\sum_{j=1}^{k_{n}}\int_{\mathcal{B}_{j}}\left[  \left(zI-\Gamma_{n}\right)  ^{-1}\left(  \Gamma_{n}-\Gamma\right)  \left(zI-\Gamma\right)  ^{-1}\right]  dz\\
&  =\mathcal{S}_{n}+\mathcal{R}_{n},
\end{align*}
where
$$
\mathcal{S}_{n}=\dfrac{1}{2\iota\pi}\sum_{j=1}^{k_{n}}\int_{\mathcal{B}_{j}}\left[  \left(  zI-\Gamma\right)  ^{-1}\left(  \Gamma_{n}-\Gamma\right)\left(  zI-\Gamma\right)  ^{-1}\right]  dz,
$$
and
\begin{equation}
\mathcal{R}_{n}=\dfrac{1}{2\iota\pi}\sum_{j=1}^{k_{n}}\int_{\mathcal{B}_{j}}\left[  \left(  zI-\Gamma\right)  ^{-1}\left(  \Gamma_{n}-\Gamma\right)\left(  zI-\Gamma\right)  ^{-1}\left(  \Gamma_{n}-\Gamma\right)  \left(zI-\Gamma_{n}\right)  ^{-1}\right]  dz. \label{residu}
\end{equation}
Result (\ref{Rn}) below will provide us with a sufficient condition for $\mathcal{R}_{n}$ to be negligible. At first, we turn to $\mathcal{S}_{n}$. We have
\begin{align*}
\mathbb{E}\left\langle \mathcal{S}_{n}\rho,X_{n+1}\right\rangle ^{2}  &  =\mathbb{E}\left(  \sum_{l,l^{\prime}=1}^{+\infty}\left\langle\mathcal{S}_{n}\rho,e_{l}\right\rangle \left\langle X_{n+1},e_{l}\right\rangle\left\langle \mathcal{S}_{n}\rho,e_{l^{\prime}}\right\rangle \left\langle X_{n+1},e_{l^{\prime}}\right\rangle \right) \\
&  =\mathbb{E}\left(  \sum_{l=1}^{+\infty}\left\langle \mathcal{S}_{n}\rho,e_{l}\right\rangle ^{2}\left\langle X_{n+1},e_{l}\right\rangle^{2}\right)  =\left(  \sum_{l=1}^{+\infty}\lambda_{l}\mathbb{E}\left\langle\mathcal{S}_{n}\rho,e_{l}\right\rangle ^{2}\right)  ,
\end{align*}
since $\mathbb{E}\left(  \left\langle X_{n+1},e_{l}\right\rangle \left\langle X_{n+1},e_{l^{\prime}}\right\rangle \right)  =0$ if $l\neq l^{\prime}$ and $X_{n+1}$ is independent from $\mathcal{S}_{n}.$ Now
$$
\mathbb{E}\left\langle \mathcal{S}_{n}\rho,e_{l}\right\rangle ^{2} =\mathbb{E}\left\langle \rho,\mathcal{S}_{n}e_{l}\right\rangle ^{2}  =\mathbb{E}\left(  \sum_{l^{\prime}=1}^{+\infty}\left\langle \rho,e_{l^{\prime}}\right\rangle \left\langle \mathcal{S}_{n}e_{l},e_{l^{\prime}}\right\rangle \right)  ^{2}.
$$
The operator $\mathcal{S}_{n}$ was explicitly computed by Dauxois \textit{et al.} (1982). More precisely
$$ \dfrac{1}{2\pi\iota}\int_{\mathcal{B}_{j}}\left[  \left(  zI-\Gamma\right)^{-1}\left(  \Gamma_{n}-\Gamma\right)  \left(  zI-\Gamma\right)  ^{-1}\right]
dz =v_{j}\left(  \Gamma_{n}-\Gamma\right)  \pi_{j}+\pi_{j}\left(  \Gamma_{n}-\Gamma\right)  v_{j},
$$
with $v_{j}=\sum_{j^{\prime}\neq j}\dfrac{1}{\lambda_{j^{\prime}}-\lambda_{j}}\pi_{j^{\prime}}$ where $\pi_{j}$ is the projector on the eigenspace associated to the $j^{th}$ eigenfunction of $\Gamma$. Hence
\begin{align}
\left\langle \mathcal{S}_{n}e_{l},e_{l^{\prime}}\right\rangle  &  =\sum_{j=1}^{k_{n}}\left[  \left\langle \left(  \Gamma_{n}-\Gamma\right)  \pi_{j}e_{l},v_{j}e_{l^{\prime}}\right\rangle +\left\langle \left(  \Gamma_{n}-\Gamma\right)  v_{j}e_{l},\pi_{j}e_{l^{\prime}}\right\rangle \right]
\nonumber\\
&  =\left\{
\begin{tabular}[c]{l}
$0$ if ($l^{\prime}\leq k_{n}$ and $l\leq k_{n}$) or if ($l^{\prime}>k_{n}$ and $l>k_{n}$),\\
$\dfrac{\left\langle \left(  \Gamma_{n}-\Gamma\right)  e_{l},e_{l^{\prime}}\right\rangle }{\lambda_{l^{\prime}}-\lambda_{l}}$ if $l^{\prime}>k_{n}$ and $l\leq k_{n}$,\\
$\dfrac{\left\langle \left(  \Gamma_{n}-\Gamma\right)  e_{l},e_{l^{\prime}}\right\rangle }{\lambda_{l}-\lambda_{l^{\prime}}}$ if $l^{\prime}\leq k_{n}$ and $l>k_{n}.$
\end{tabular}
\ \ \ \ \ \ \ \right.  \label{triple}
\end{align}
Finally, if we take for instance $l\leq k_{n}$
\begin{align*}
\mathbb{E}\left\langle \mathcal{S}_{n}\rho,e_{l}\right\rangle ^{2}  &  =\mathbb{E}\left(  \sum_{l^{\prime}\geq k_{n}+1}^{+\infty}\left\langle\rho,e_{l^{\prime}}\right\rangle \dfrac{\left\langle \left(  \Gamma_{n}-\Gamma\right)  e_{l},e_{l^{\prime}}\right\rangle }{\lambda_{l^{\prime}}-\lambda_{l}}\right)  ^{2}\\
&  =\mathbb{E}\left(  \dfrac{1}{n}\sum_{j=1}^{n}\sum_{l^{\prime}\geq k_{n}+1}^{+\infty}\left\langle \rho,e_{l^{\prime}}\right\rangle \dfrac{\left\langle (X_{j}\otimes X_j-\Gamma)e_{l},e_{l^{\prime}}\right\rangle }{\lambda_{l^{\prime}}-\lambda_{l}}\right)  ^{2}\\
&  =\mathbb{E}\left(  \dfrac{1}{n}\sum_{j=1}^{n}Z_{j,l,n}^{\ast}\right)  ^{2},
\end{align*}
where
$$
Z_{j,l,n}^{\ast}=\sum_{l^{\prime}\geq k_{n}+1}^{+\infty}\left\langle\rho,e_{l^{\prime}}\right\rangle \dfrac{\left\langle \left(  X_{j}\otimes X_{j}-\Gamma\right)  e_{l},e_{l^{\prime}}\right\rangle }{\lambda_{l^{\prime}}-\lambda_{l}},
$$
and the $\left(  Z_{j,l,n}^{\ast}\right)  _{j\geq1}$ are centered and uncorrelated random variables. Hence
$$
\mathbb{E}\left(  \dfrac{1}{n}\sum_{j=1}^{n}Z_{j,l,n}^{\ast}\right)  ^{2}  =\dfrac{1}{n}\mathbb{E}\left(  \sum_{l^{\prime}\geq k_{n}+1}^{+\infty}\left\langle \rho,e_{l^{\prime}}\right\rangle \dfrac{\left\langle \left(X_{1},e_{l}\right)  \right\rangle \left\langle \left(  X_{1},e_{l^{\prime}}\right)  \right\rangle }{\lambda_{l^{\prime}}-\lambda_{l}}\right)  ^{2}.
$$
Since $l\leq k_{n}<l^{\prime}$, by using the Karhunen-Loève expansion of $X_{1}$, we get
$$
\sum_{l^{\prime}\geq k_{n}+1}^{+\infty}\left\langle \rho,e_{l^{\prime}}\right\rangle \dfrac{\left\langle \left(  X_{1},e_{l}\right)  \right\rangle\left\langle \left(  X_{1},e_{l^{\prime}}\right)  \right\rangle }{\lambda_{l^{\prime}}-\lambda_{l}}=\sum_{l^{\prime}\geq k_{n}+1}^{+\infty}\left\langle \rho,e_{l^{\prime}}\right\rangle \dfrac{\sqrt{\lambda_{l}\lambda_{l^{\prime}}}\xi_{l}\xi_{l^{\prime}}}{\lambda_{l^{\prime}}-\lambda_{l}}.
$$
and then
$$
\mathbb{E}\left(  \sum_{l^{\prime}\geq k_{n}+1}^{+\infty}\left\langle\rho,e_{l^{\prime}}\right\rangle \dfrac{\left\langle \left(  X_{1},e_{l}\right)  \right\rangle \left\langle \left(  X_{1},e_{l^{\prime}}\right)\right\rangle }{\lambda_{l^{\prime}}-\lambda_{l}}\right)  ^{2}=\sum_{l^{\prime},m\geq k_{n}+1}^{+\infty}\left\langle \rho,e_{l^{\prime}}\right\rangle \left\langle \rho,e_{m}\right\rangle \dfrac{\sqrt{\lambda_{l}^{2}\lambda_{l}\lambda_{m}}\mathbb{E}\left(  \xi_{l}^{2}\xi_{l^{\prime}}\xi_{m}\right)  }{\left(  \lambda_{l^{\prime}}-\lambda_{l}\right)  \left(\lambda_{m}-\lambda_{l}\right)  }.
$$
By applying twice Cauchy-Schwarz inequality to the $\xi_{k}$'s and under assumption $(A.3)$, we get
\begin{align*}
\mathbb{E}\left(  \xi_{l}^{2}\xi_{l^{\prime}}\xi_{m}\right)   &  \leq \sqrt{\mathbb{E}\left(  \xi_{l}^{4}\right)  }\sqrt{\mathbb{E}\left(\xi_{l^{\prime}}^{2}\xi_{m}^{2}\right)  }\\
&  \leq\sqrt{M}\sqrt{M}.
\end{align*}
Summing up what we made above we get
$$
\mathbb{E}\left(  \sum_{l^{\prime}\geq k_{n}+1}^{+\infty}\left\langle\rho,e_{l^{\prime}}\right\rangle \dfrac{\left\langle \left(  X_{1},e_{l}\right)  \right\rangle \left\langle \left(  X_{1},e_{l^{\prime}}\right)\right\rangle }{\lambda_{l^{\prime}}-\lambda_{l}}\right)  ^{2}\leq M\left(\sum_{l^{\prime}\geq k_{n}+1}^{+\infty}\left\langle \rho,e_{l^{\prime}}\right\rangle \dfrac{\sqrt{\lambda_{l}\lambda_{l^{\prime}}}}{\lambda_{l^{\prime}}-\lambda_{l}}\right)  ^{2}.
$$
Remember that we had fixed $l\leq k_{n}$. Now, if we take $l>k_{n}$ similar calculations lead to
$$
\mathbb{E}\left(  \sum_{l^{\prime}=1}^{k_{n}}\left\langle \rho,e_{l^{\prime}}\right\rangle \dfrac{\left\langle \left(  X_{1},e_{l}\right)  \right\rangle\left\langle \left(  X_{1},e_{l^{\prime}}\right)  \right\rangle }{\lambda_{l^{\prime}}-\lambda_{l}}\right)  ^{2}\leq M\left(  \sum_{l^{\prime}=1}^{k_{n}}\left\langle \rho,e_{l^{\prime}}\right\rangle \dfrac{\sqrt{\lambda_{l}\lambda_{l^{\prime}}}}{\lambda_{l^{\prime}}-\lambda_{l}}\right)^{2}.
$$
At last
\begin{align}
\dfrac{n}{k_{n}}\mathbb{E}\left\langle \mathcal{S}_{n}\rho,X_{n+1}\right\rangle ^{2}  &  \leq\dfrac{M}{k_{n}}\sum_{l=1}^{k_{n}}\lambda_{l}\left(  \sum_{l^{\prime}\geq k_{n}+1}^{+\infty}\left\langle \rho,e_{l^{\prime}}\right\rangle \dfrac{\sqrt{\lambda_{l}\lambda_{l^{\prime}}}}{\lambda_{l^{\prime}}-\lambda_{l}}\right)  ^{2}\label{nosouci1}\\
&  +\dfrac{M}{k_{n}}\sum_{l>k_{n}}\lambda_{l}\left(  \sum_{l^{\prime}\geq1}^{k_{n}}\left\langle \rho,e_{l^{\prime}}\right\rangle \dfrac{\sqrt{\lambda_{l}\lambda_{l^{\prime}}}}{\lambda_{l^{\prime}}-\lambda_{l}}\right)^{2}. \label{nosouci2}
\end{align}
We apply Lemma \ref{trick1} first to bound (\ref{nosouci1})
\begin{align*}
\dfrac{M}{k_{n}}\sum_{l=1}^{k_{n}}\lambda_{l}\left(  \sum_{l^{\prime}\geq k_{n}+1}^{+\infty}\left\langle \rho,e_{l^{\prime}}\right\rangle \dfrac{\sqrt{\lambda_{l}\lambda_{l^{\prime}}}}{\lambda_{l}-\lambda_{l^{\prime}}}\right)  ^{2} &  \leq\dfrac{M}{k_{n}}\sum_{l=1}^{k_{n}}\lambda_{l}\left(  \sum_{l^{\prime}\geq k_{n}+1}^{+\infty}\left\langle \rho,e_{l^{\prime}}\right\rangle\sqrt{\dfrac{\lambda_{l^{\prime}}}{\lambda_{l}}}\dfrac{1}{1-\dfrac{l}{l^{\prime}}}\right)  ^{2}\\
&  \leq\dfrac{M}{k_{n}}\sum_{l=1}^{k_{n}}\left(  \sum_{l^{\prime}\geq k_{n}+1}^{+\infty}\left\langle \rho,e_{l^{\prime}}\right\rangle \sqrt
{\lambda_{l^{\prime}}}\dfrac{1}{1-\dfrac{l}{l^{\prime}}}\right)  ^{2}.
\end{align*}
Now we set $h_{n}=\left[  \sqrt{\dfrac{k_{n}}{\log k_{n}}}\right]  $ where $\left[  u\right]  ,$ $u\in\mathbb{R}$, denotes the largest integer smaller
than $u.$ Note that the last inequality in the display above may be split as follows
\begin{align}
\dfrac{M}{k_{n}}\sum_{l=1}^{k_{n}}\lambda_{l}\left(  \sum_{l^{\prime}\geq k_{n}+1}^{+\infty}\left\langle \rho,e_{l^{\prime}}\right\rangle \dfrac{\sqrt{\lambda_{l}\lambda_{l^{\prime}}}}{\lambda_{l}-\lambda_{l^{\prime}}}\right)  ^{2} &  \leq\dfrac{2M}{k_{n}}\sum_{l=1}^{k_{n}}\left(  \sum_{l^{\prime}\geq k_{n}+1}^{k_{n}+h_{n}}\left\vert \left\langle \rho,e_{l^{\prime}}\right\rangle\right\vert \sqrt{\lambda_{l^{\prime}}}\dfrac{1}{1-\dfrac{l}{l^{\prime}}}\right)  ^{2}\label{qqq}\\
&  +\dfrac{2M}{k_{n}}\sum_{l=1}^{k_{n}}\left(  \sum_{l^{\prime}\geq k_{n}+h_{n}}^{+\infty}\left\vert \left\langle \rho,e_{l^{\prime}}\right\rangle\right\vert \sqrt{\lambda_{l^{\prime}}}\dfrac{1}{1-\dfrac{l}{l^{\prime}}}\right)  ^{2}.\nonumber
\end{align}
Dealing with the second term we get for $l^{\prime}\geq k_{n}+h_{n}$
$$
1-\dfrac{l}{l^{\prime}}\geq1-\dfrac{k_{n}}{k_{n}+h_{n}}=\dfrac{h_{n}}%
{k_{n}+h_{n}},
$$
and hence
\begin{align*}
\sum_{l^{\prime}\geq k_{n}+h_{n}}^{+\infty}\left\vert \left\langle\rho,e_{l^{\prime}}\right\rangle \right\vert \sqrt{\lambda_{l^{\prime}}}\dfrac{1}{1-\dfrac{l}{l^{\prime}}}  &  \leq\sum_{l^{\prime}\geq k_{n}+h_{n}}^{+\infty}\left\vert \left\langle \rho,e_{l^{\prime}}\right\rangle\right\vert \sqrt{\lambda_{l^{\prime}}}\left(  1+\dfrac{k_{n}}{h_{n}}\right)
\\
&  \leq\sum_{l^{\prime}\geq k_{n}+h_{n}}^{+\infty}\left\vert \left\langle\rho,e_{l^{\prime}}\right\rangle \right\vert \sqrt{\lambda_{l^{\prime}}}\left(  1+\sqrt{k_{n}\log k_{n}}\right)  .
\end{align*}
Now obviously
$$
\sup_{l^{\prime}\geq k_{n}+h_{n}}\sqrt{\lambda_{l^{\prime}}}\left(1+\sqrt{k_{n}\log k_{n}}\right)  \leq K,
$$
since $\sqrt{\lambda_{l^{\prime}}l^{\prime}\log l^{\prime}}\rightarrow0$ from which we deduce that
$$
\dfrac{2M}{k_{n}}\sum_{l=1}^{k_{n}}\left(  \sum_{l^{\prime}\geq k_{n}+h_{n}}^{+\infty}\left\langle \rho,e_{l^{\prime}}\right\rangle \sqrt{\lambda_{l^{\prime}}}\dfrac{1}{1-\dfrac{l}{l^{\prime}}}\right)  ^{2}\leq\dfrac{2MK}{k_{n}}\sum_{l=1}^{k_{n}}\left(  \sum_{l^{\prime}\geq k_{n}+h_{n}}^{+\infty}\left\vert \left\langle \rho,e_{l^{\prime}}\right\rangle\right\vert \right)  ^{2}.
$$
When assumption $(A.3)$ holds, Cesaro's mean Theorem ensures that the term on the left above tends to zero. We turn to the first term in equation (\ref{qqq})
\begin{align*}
\left(  \sum_{l^{\prime}\geq k_{n}+1}^{k_{n}+h_{n}}\left\vert \left\langle\rho,e_{l^{\prime}}\right\rangle \right\vert \sqrt{\lambda_{l^{\prime}}}\dfrac{1}{1-\dfrac{l}{l^{\prime}}}\right)  ^{2}  &  \leq h_{n}^{2}\max_{\substack{k_{n}+1\leq l^{\prime}\leq k_{n}+h_{n},\\1\leq l\leq k_{n}}}\left\{  \left\vert \left\langle \rho,e_{l^{\prime}}\right\rangle\right\vert \sqrt{\lambda_{l^{\prime}}}\dfrac{1}{1-\dfrac{l}{l^{\prime}}}\right\}  ^{2}\\
&  \leq\dfrac{k_{n}}{\log k_{n}}\lambda_{k_{n}}k_{n}^{2}\max_{k_{n}+1\leq l^{\prime}\leq k_{n}+h_{n},}\left(  \left\vert \left\langle \rho,e_{l^{\prime}}\right\rangle \right\vert ^{2}\right)  .
\end{align*}
Now $\lambda_{k_{n}}k_{n}$ as well as $k_{n}\max_{k_{n}+1\leq l^{\prime}\leq k_{n}+h_{n},}\left(  \left\vert \left\langle \rho,e_{l^{\prime}}\right\rangle \right\vert \right)  $ tend to zero when assumption $(A.3)$ holds. We get once more
$$
\dfrac{2M}{k_{n}}\sum_{l=1}^{k_{n}}\left(  \sum_{l^{\prime}\geq k_{n}+1}^{k_{n}+h_{n}}\left\vert \left\langle \rho,e_{l^{\prime}}\right\rangle\right\vert \sqrt{\lambda_{l^{\prime}}}\dfrac{1}{1-\dfrac{l}{l^{\prime}}}\right)  ^{2}\rightarrow0.
$$
A similar truncating technique would prove that the term in (\ref{nosouci2}) also tends to zero as $n$ goes to infinity which leads to
\begin{equation}
\dfrac{n}{k_{n}}\mathbb{E}\left\langle \mathcal{S}_{n}\rho,X_{n+1}\right\rangle ^{2}\rightarrow0. \label{alti}
\end{equation}
In order to finish the proof of the Proposition we must deal with the term
introduced in (\ref{residu}). We have the following result
\begin{equation}\label{Rn}
\sqrt{\dfrac{n}{k_{n}}}\left\vert \left\langle \mathcal{R}_{n}\rho,X_{n+1}\right\rangle \right\vert =O_{\mathbb{P}}\left(  \dfrac{1}{\sqrt{n}}k_{n}^{5/2}\left(  \log k_{n}\right)  ^{2}\right),
\end{equation}
when $\dfrac{k_{n}^{2}\log k_{n}}{\sqrt{n}}\rightarrow0$. Indeed, consider
$$
T_{j,n}=\int_{\mathcal{B}_{j}}\left[  \left(  zI-\Gamma\right)  ^{-1}\left(\Gamma_{n}-\Gamma\right)  \left(  zI-\Gamma\right)  ^{-1}\left(  \Gamma_{n}-\Gamma\right)  \left(  zI-\Gamma_{n}\right)  ^{-1}\right]  dz.
$$
Then setting
$$
G_{n}\left(  z\right)  =\left(  zI-\Gamma\right)  ^{-1/2}\left(\Gamma_{n}-\Gamma\right)  \left(  zI-\Gamma\right)  ^{-1/2},
$$
we have
\begin{align}
&  \left\vert \left\langle T_{j,n}\rho,X_{n+1}\right\rangle \right\vert\label{ketchup}\\
&  =\left\vert \int_{\mathcal{B}_{j}}\left\langle \left(  zI-\Gamma\right)^{-1/2}\left(  \Gamma_{n}-\Gamma\right)  \left(  zI-\Gamma\right)^{-1}\left(  \Gamma_{n}-\Gamma\right)  \left(  zI-\Gamma_{n}\right)  ^{-1}\rho,\left(  zI-\Gamma\right)  ^{-1/2}X_{n+1}\right\rangle dz\right\vert\nonumber\\
&  \leq\int_{\mathcal{B}_{j}}\left\vert \left\langle G_{n}^{2}\left(z\right)  \left(  zI-\Gamma\right)  ^{1/2}\left(  zI-\Gamma_{n}\right)^{-1}\left(  zI-\Gamma\right)  ^{1/2}\left(  zI-\Gamma\right)  ^{-1/2}\rho,\left(  zI-\Gamma\right)  ^{-1/2}X_{n+1}\right\rangle \right\vert dz\nonumber\\
&  \leq\int_{\mathcal{B}_{j}}\left\Vert G_{n}\left(  z\right)  \right\Vert_{\infty}^{2}\left\Vert \left(  zI-\Gamma\right)  ^{1/2}\left(  zI-\Gamma_{n}\right)  ^{-1}\left(  zI-\Gamma\right)  ^{1/2}\right\Vert _{\infty}\left\Vert \left(  zI-\Gamma\right)  ^{-1/2}X_{n+1}\right\Vert \left\Vert\left(  zI-\Gamma\right)  ^{-1/2}\rho\right\Vert dz.\nonumber
\end{align}
Following Lemma \ref{you}, the random variable $\left\Vert \left(zI-\Gamma\right)  ^{1/2}\left(  zI-\Gamma_{n}\right)  ^{-1}\left(zI-\Gamma\right)  ^{1/2}\right\Vert _{\infty}$ is decomposed in two terms
$$
\left\Vert \left(  zI-\Gamma\right)  ^{1/2}\left(  zI-\Gamma_{n}\right)^{-1}\left(  zI-\Gamma\right)  ^{1/2}\right\Vert _{\infty}\left(1\!\!1_{\mathcal{E}_{j}\left(  z\right)  }+1\!\!1_{\mathcal{E}_{j}^{c}\left(z\right)  }\right)  .
$$
On the one hand when $\mathcal{E}_{j}\left(  z\right)  $ holds it was proved in Lemma \ref{you} that
\begin{equation}
\left\Vert \left(  zI-\Gamma\right)  \right\Vert _{\infty}\left\Vert \left(zI-\Gamma_{n}\right)  ^{-1}\right\Vert _{\infty}\leq C. \label{gio}
\end{equation}
On the other hand when $\mathcal{E}_{j}^{c}\left(  z\right)  $ holds we may write for all $\eta>0$ thanks to bound (\ref{anytime})
$$
\mathbb{P}\left(  \left\vert \left\langle T_{j,n}\rho,X_{n+1}\right\rangle\right\vert 1\!\!1_{\mathcal{E}_{j}^{c}\left(  z\right)  }>\eta\right)\leq\mathbb{P}\left(  \mathcal{E}_{j}^{c}\left(  z\right)  \right)  \leq \dfrac{M}{\sqrt{n}}\left(  j\log j\right),
$$
which entails that
\begin{align*}
&  \mathbb{P}\left(  \sum_{j=1}^{k_{n}}\left\vert \left\langle T_{j,n}\rho,X_{n+1}\right\rangle \right\vert 1\!\!1_{\mathcal{E}_{j}^{c}\left(z\right)  }>\eta\right) \\
&  \leq M\sum_{j=1}^{k_{n}}\dfrac{1}{\sqrt{n}}\left(  j\log j\right)\leq\dfrac{k_{n}^{2}\log k_{n}}{\sqrt{n}}\rightarrow0.
\end{align*}
\newline Consequently we can deal with all $T_{j,n}$ as if the event $\mathcal{E}_{j}\left(  z\right)  $ -hence the bound (\ref{gio})- holds almost surely. We take expectation and note that $G_{n}\left(  z\right)  $ and $X_{n+1}$ are independent
$$\mathbb{E}\left\vert \left\langle T_{j,n}\rho,X_{n+1}\right\rangle\right\vert \leq C\int_{\mathcal{B}_{j}}\mathbb{E}\left\Vert G_{n}\left(  z\right)\right\Vert _{\infty}^{2}\mathbb{E}\left\Vert \left(  zI-\Gamma\right)^{-1/2}X_{n+1}\right\Vert \left\Vert \left(  zI-\Gamma\right)  ^{-1/2}\rho\right\Vert dz.
$$
By Lemma \ref{thal} we have
\begin{align*}
\mathbb{E}\left\vert \left\langle T_{j,n}\rho,X_{n+1}\right\rangle\right\vert  &  \leq\dfrac{C}{n}\mathrm{diam}\left(  \mathcal{B}_{j}\right)\cdot\left(  j\log j\right)  ^{5/2}\sup_{z\in\mathcal{B}_{j}}\left\Vert\left(  zI-\Gamma\right)  ^{-1/2}\rho\right\Vert \\
&  \leq C\sqrt{\delta_{j}}\cdot\left(  j\log j\right)  ^{5/2}\left\Vert\rho\right\Vert \leq C\left(  j\log j\right)  ^{2},
\end{align*}
since $\delta_{j}\leq C\left(  j\log j\right)  ^{-1}$ at least for a sufficiently large $j.$ Finally summing over all the $j$'s from $1$ to $k_{n}$ leads to
$$
\sqrt{\dfrac{n}{k_{n}}}\mathbb{E}\left\vert \left\langle \mathcal{R}_{n}\rho,X_{n+1}\right\rangle \right\vert \leq C\dfrac{1}{\sqrt{nk_{n}}}\sum_{j=1}^{k_{n}}\left(  j\log j\right)  ^{2}\leq\dfrac{C}{\sqrt{n}}k_{n}^{5/2}\left(  \log k_{n}\right)  ^{2},
$$
which proves (\ref{Rn}) and achieves the proof of the proposition.
\end{proof}

The methods used to prove the next Proposition are close to those developed above.

\begin{proposition}
\label{dd}If $\dfrac{1}{\sqrt{n}}k_{n}^{5/2}\left(  \log k_{n}\right)^{2}\rightarrow0$, then
$$
\sqrt{\dfrac{n}{k_{n}}}\left\vert \left\langle \left(  \Gamma_{n}^{\dagger}-\Gamma^{\dagger}\right)  U_{n},X_{n+1}\right\rangle \right\vert\overset{\mathbb{P}}{\rightarrow}0.
$$
Besides if $x$ is a fixed vector in $H$ such that $\sup_{p}\dfrac{\left\vert\left\langle x,e_{p}\right\rangle \right\vert ^{2}}{\lambda_{p}}<+\infty$ and $\dfrac{k_{n}^{3}\left(  \log k_{n}\right)  ^{2}}{t_{n,x}\sqrt{n}}\rightarrow0,$
$$
\dfrac{\sqrt{n}}{t_{n,x}}\left\vert \left\langle \left(  \Gamma_{n}^{\dagger}-\Gamma^{\dagger}\right)  U_{n},x\right\rangle \right\vert \overset{\mathbb{P}}{\rightarrow}0.
$$
\end{proposition}

\begin{proof}
Once again we develop the expression above by means of complex integrals for operator-valued analytic functions. Hence
\begin{align*}
\Gamma_{n}^{\dagger}-\Gamma^{\dagger}  &  =\dfrac{1}{2\pi\iota}\int_{\mathcal{C}_{n}}\widetilde{f}_{n}\left(z\right)  \left[  \left(  zI-\Gamma_{n}\right)  ^{-1}\left(  \Gamma-\Gamma_{n}\right)  \left(  zI-\Gamma\right)  ^{-1}\right]  dz\\
&  =\sum_{j=1}^{k_{n}}\dfrac{1}{2\pi\iota}\int_{\mathcal{B}_{j}}\widetilde{f}_{n}\left(  z\right)  \left[  \left(  zI-\Gamma_{n}\right)  ^{-1}\left(\Gamma-\Gamma_{n}\right)  \left(  zI-\Gamma\right)  ^{-1}\right]  dz,
\end{align*}
and
$$
\left\vert \left\langle \left(  \Gamma_{n}^{\dagger}-\Gamma^{\dagger}\right)U_{n},X_{n+1}\right\rangle \right\vert \leq C\sum_{j=1}^{k_{n}}H_{j,n},
$$
where
$$
H_{j,n}=\int_{\mathcal{B}_{j}}\left\vert \widetilde{f}_{n}\left(  z\right)\left\langle \left(  zI-\Gamma\right)  ^{1/2}\left(  zI-\Gamma_{n}\right)^{-1}\left(  zI-\Gamma\right)  ^{1/2}G_{n}\left(  z\right)  \left(zI-\Gamma\right)  ^{-1/2}U_{n},\left(  zI-\Gamma\right)  ^{-1/2}X_{n+1}\right\rangle \right\vert dz.
$$
We copy verbatim the arguments used to bound (\ref{ketchup}) : first of all we reintroduce the operator $G_{n}\left(  z\right)  $ below and
$$
\left(  zI-\Gamma\right)  ^{1/2}\left(  zI-\Gamma_{n}\right)  ^{-1}\left(zI-\Gamma\right)  ^{1/2},
$$
remains almost surely bounded by a constant which does not depend on $n$ or $j$ plus a negligible term as was proved just below (\ref{gio}). Hence
$$
H_{j,n}\leq C\int_{\mathcal{B}_{j}}\left\vert \widetilde{f}_{n}\left(z\right)  \right\vert \left\Vert G_{n}\left(  z\right)  \right\Vert \left\Vert\left(  zI-\Gamma\right)  ^{-1/2}U_{n}\right\Vert \left\Vert \left(zI-\Gamma\right)  ^{-1/2}X_{n+1}\right\Vert dz.
$$
We take expectation
\begin{align*}
\mathbb{E}H_{j,n}  &  \leq C\int_{\mathcal{B}_{j}}\left\vert \widetilde{f}_{n}\left(  z\right)  \right\vert \mathbb{E}\left(  \left\Vert G_{n}\left(z\right)  \right\Vert \left\Vert \left(  zI-\Gamma\right)  ^{-1/2}U_{n}\right\Vert \right)  \mathbb{E}\left\Vert \left(  zI-\Gamma\right)
^{-1/2}X_{n+1}\right\Vert dz\\
&  \leq C\mathrm{diam}\left(  \mathcal{B}_{j}\right)  \sup_{z\in\mathcal{B}_{j}}\left(  \left\vert \widetilde{f}_{n}\left(  z\right)\right\vert \mathbb{E}\left\Vert \left(  zI-\Gamma\right)  ^{-1/2}X_{n+1}\right\Vert \sqrt{\mathbb{E}\left\Vert G_{n}\left(  z\right)\right\Vert ^{2}}\sqrt{\mathbb{E}\left\Vert \left(  zI-\Gamma\right)^{-1/2}U_{n}\right\Vert ^{2}}\right),
\end{align*}
where Cauchy-Schwarz inequality was applied. Now invoking Lemma \ref{thal} yields
$$
\mathbb{E}H_{j,n}\leq C\dfrac{\mathrm{diam}\left(  \mathcal{B}_{j}\right)}{\sqrt{n}}\left(  j\log j\right)  ^{3/2}\sup_{z\in\mathcal{B}_{j}}\left(\left\vert \widetilde{f}_{n}\left(  z\right)  \right\vert \sqrt{\mathbb{E}\left\Vert \left(  zI-\Gamma\right)  ^{-1/2}U_{n}\right\Vert ^{2}}\right).
$$
Obviously
\begin{align*}
\mathbb{E}\left\Vert \left(  zI-\Gamma\right)  ^{-1/2}U_{n}\right\Vert ^{2}&  =\dfrac{\sigma_{\varepsilon}^{2}}{n}\mathbb{E}\left\Vert \left(zI-\Gamma\right)  ^{-1/2}X_{1}\right\Vert ^{2}\\
&  =\dfrac{\sigma_{\varepsilon}^{2}}{n}\sum_{l=1}^{+\infty}\dfrac{\lambda_{l}}{\left\vert z-\lambda_{l}\right\vert },
\end{align*}
hence
$$
\sup_{z\in\mathcal{B}_{j}}\left(  \sqrt{\mathbb{E}\left\Vert \left(zI-\Gamma\right)  ^{-1/2}U_{n}\right\Vert ^{2}}\right)  \leq\dfrac{1}{\sqrt{n}}\left(  j\log j\right)  ^{1/2}.
$$
At last
$$
\mathbb{E}H_{j,n}\leq C\dfrac{\delta_{j}}{\lambda_{j}n}\left(  j\log j\right)^{2}\leq\dfrac{C}{n}\left(  j\log j\right)  ^{2},
$$
and
$$
\mathbb{E}\left\vert \left\langle \left(  \Gamma_{n}^{\dagger}-\Gamma^{\dagger}\right)  U_{n},X_{n+1}\right\rangle \right\vert \leq\dfrac{C}{n}k_{n}^{3}\left(  \log k_{n}\right)  ^{2},
$$
which proves the first part of the Proposition. Replacing $X_{n+1}$ with a fixed $x$ in $H,$ means replacing $\mathbb{E}\left\Vert \left(zI-\Gamma\right)  ^{-1/2}X_{n+1}\right\Vert $ with
$$
\left\Vert \left(  zI-\Gamma\right)  ^{-1/2}x\right\Vert \leq\sqrt{\sum_{p=1}^{+\infty}\dfrac{\left\langle x,e_{p}\right\rangle ^{2}}{\left\vert z-\lambda_{p}\right\vert }}\leq\sqrt{\sup_{p}\dfrac{\left\vert \left\langle x,e_{p}\right\rangle \right\vert }{\lambda_{p}}}\sqrt{\sum_{p=1}^{+\infty}\dfrac{\lambda_{p}}{\left\vert z-\lambda_{p}\right\vert }},
$$
and the derivation of the second part of the Proposition stems from the first part.
\end{proof}


\subsection{Weakly convergent terms}


This subsection is quite short but was separated from the others for the sake of clarity and in order to give a logical structure to the proofs.

\begin{lemma}
\label{Rn0}We have
$$
\sqrt{\frac{n}{t_{n,x}}}\left\langle R_{n},x\right\rangle \overset{w}{\rightarrow}N\left(  0,\sigma_{\varepsilon}^{2}\right)  , x\in H,
$$
and
$$
\sqrt{\frac{n}{s_{n}}}\left\langle R_{n},X_{n+1}\right\rangle \overset{w}{\rightarrow}N\left(  0,\sigma_{\varepsilon}^{2}\right).
$$
\end{lemma}

\begin{proof} We have
$$
\left\langle R_{n},x\right\rangle =\left\langle \Gamma^{\dagger}U_{n},x\right\rangle =\dfrac{1}{n}\sum_{i=1}^{n}\left\langle \Gamma^{\dagger}X_{i},x\right\rangle \varepsilon_{i},
$$
which is an array - $\Gamma^{\dagger}$ implicitly depends on $n$ - of independent real r.v. The Central Limit Theorem holds for this sequence and leads to the first announced result. We turn to the second display
\begin{align*}
\left\langle R_{n},X_{n+1}\right\rangle  &  =\left\langle \Gamma^{\dagger}U_{n},X_{n+1}\right\rangle \\
&  =\dfrac{1}{n}\sum_{i=1}^{n}\left\langle \Gamma^{\dagger}X_{i},X_{n+1}\right\rangle \varepsilon_{i}=\sum_{i=1}^{n}Z_{i,n}.
\end{align*}
Denoting $\mathcal{F}_{i}$ the $\sigma$-algebra generated by $\left(X_{1},\varepsilon_{1},...,X_{i},\varepsilon_{i}\right)$, we see that$Z_{i,n}$ is a martingale difference sequence w.r.t. $\mathcal{F}_{i}$. Also note that
$$
\mathbb{E}\left(  Z_{i,n}^{2}|\mathcal{F}_{i}\right)  =\dfrac{\varepsilon_{i}^{2}}{n^{2}}\left\Vert \Gamma^{1/2}\Gamma^{\dagger}X_{i}\right\Vert ^{2},
$$
and that
\begin{align*}
\mathbb{E}\left[  \varepsilon_{i}^{2}\left\Vert \Gamma^{1/2}\Gamma^{\dagger}X_{i}\right\Vert ^{2}\right] & =\mathbb{E}\left[  \left\Vert \Gamma^{1/2}\Gamma^{\dagger}X_{i}\right\Vert^{2}\mathbb{E}\left(  \varepsilon_{i}^{2}|X_{i}\right)  \right] \\
&  =\sigma_{\varepsilon}^{2}\mathbb{E}\left\Vert \Gamma^{1/2}\Gamma^{\dagger}X_{i}\right\Vert ^{2}\\
&  =\sigma_{\varepsilon}^{2}\sum_{j=1}^{k_{n}}\left[  \lambda_{j}f_{n}\left(\lambda_{j}\right)  \right]  ^{2}=\sigma_{\varepsilon}^{2}s_{n}^{2}.
\end{align*}
Applying the Central Limit Theorem for real valued martingale difference arrays (see e.g Mc Leish, 1974) we get the second result.
\end{proof}


\subsection{Proofs of the main results}


The careful reader has noted that within the preceding steps of the proofs $s_{n}^{2}$ was replaced with $k_{n}$ in the normalizing sequence. Very simple computations prove that under (H.2) and if $k_{n}/\sqrt{n}$ tends to zero this permutation is possible (it is enough to prove that $k_{n}\geq Cs_{n}^{2}$ for some constant $C$).

\textbf{Proof of Theorem \ref{best}}

The proof of Theorem \ref{best} stems from the decomposition (\ref{decomp}), Lemma \ref{Tn}, Proposition \ref{ks}, Proposition \ref{dd} and from Lemma \ref{Rn0}.

\textbf{Proof of Corollary \ref{bias}.}

The proof of the Corollary is a straightforward consequence of Lemma \ref{Ln} when choosing $n=k_{n}^{6}$.

\textbf{Proof of Corollary \ref{Cor2}}

It suffices to prove that $\dfrac{\left\vert \widehat{s}_{n}^{2}-s_{n}^{2}\right\vert }{s_{n}^{2}}\overset{\mathbb{P}}{\rightarrow}0,$ or equivalently that
$$
\dfrac{\sum_{j=1}^{k_{n}}\left\vert \lambda_{j}f_{n}\left(  \lambda
_{j}\right)  -\widehat{\lambda}_{j}f_{n}\left(  \widehat{\lambda}_{j}\right)
\right\vert \left(  \lambda_{j}f_{n}\left(  \lambda_{j}\right)  +\widehat
{\lambda}_{j}f_{n}\left(  \widehat{\lambda}_{j}\right)  \right)  }{\sum
_{j=1}^{k_{n}}\left[  \lambda_{j}f_{n}\left(  \lambda_{j}\right)  \right]
^{2}}\overset{\mathbb{P}}{\rightarrow}0.
$$
Clearly since $\sup_{j\in\mathbb{N}}\left\vert \widehat{\lambda}_{j}-\lambda_{j}\right\vert =O_{P}\left(  1/\sqrt{n}\right)  $ and $xf_{n}(x)$ is bounded for $x>c_{n}$ it is enough to get
\begin{equation}
\dfrac{\sum_{j=1}^{k_{n}}\left\vert \lambda_{j}f_{n}\left(  \lambda_{j}\right)  -\widehat{\lambda}_{j}f_{n}\left(  \widehat{\lambda}_{j}\right)\right\vert }{\sum_{j=1}^{k_{n}}\left[  \lambda_{j}f_{n}\left(  \lambda_{j}\right)  \right]  ^{2}}\overset{\mathbb{P}}{\rightarrow}0 \label{del}.
\end{equation}
But by assumption $(H.3)$
\begin{align*}
\sum_{j=1}^{k_{n}}\left\vert \lambda_{j}f_{n}\left(  \lambda_{j}\right)-\widehat{\lambda}_{j}f_{n}\left(  \widehat{\lambda}_{j}\right)  \right\vert &  \leq\sum_{j=1}^{k_{n}}\left\vert \lambda_{j}f_{n}\left(  \lambda_{j}\right)  -1\right\vert +\sum_{j=1}^{k_{n}}\left\vert 1-\widehat{\lambda}_{j}f_{n}\left(  \widehat{\lambda}_{j}\right)  \right\vert \\
&  =o_{P}\left(  k_{n}/\sqrt{n}\right),
\end{align*}
and $k_{n}/\sqrt{n}\rightarrow0.$ \endproof

\vspace{0.5cm}
\textbf{Proof of Theorem \ref{TH0}}

Like Theorem \ref{best}, the proof of Theorem \ref{TH0} stems from (\ref{decomp}), Proposition \ref{dd} and from Lemma \ref{Rn0}.


\textbf{Proof of Corollary \ref{Cor1}}

We have to prove that
$$
\dfrac{\widehat{t}_{n,x}^{2}-t_{n,x}^{2}}{t_{n,x}^{2}}\ =\ \dfrac{\sum_{j=1}^{k_{n}}\widehat{\lambda}_{j}\left[  f_{n}\left(  \widehat{\lambda}_{j}\right)  \right]  ^{2}\left\langle x,\widehat{e}_{j}\right\rangle^{2}-\lambda_{j}\left[  f_{n}\left(  \lambda_{j}\right)  \right]^{2}\left\langle x,e_{j}\right\rangle ^{2}}{\sum_{j=1}^{k_{n}}\lambda_{j}\left[  f_{n}\left(  \lambda_{j}\right)  \right]  ^{2}\left\langle x,e_{j}\right\rangle ^{2}}\ \overset{\mathbb{P}}{\rightarrow}\ 0.
$$
We split the expression into two terms
\begin{align*}
w_{n1}  &  =\dfrac{\sum_{j=1}^{k_{n}}\left(  \widehat{\lambda}_{j}\left[f_{n}\left(  \widehat{\lambda}_{j}\right)  ^{2}\right]  -\lambda_{j}\left[f_{n}\left(  \lambda_{j}\right)  \right]  ^{2}\right)  \left\langle x,\widehat{e}_{j}\right\rangle ^{2}}{\sum_{j=1}^{k_{n}}\lambda_{j}\left[f_{n}\left(  \lambda_{j}\right)  \right]  ^{2}\left\langle x,e_{j}\right\rangle ^{2}},\\
w_{n2}  &  =\dfrac{\sum_{j=1}^{k_{n}}\lambda_{j}\left[  f_{n}\left(\lambda_{j}\right)  \right]  ^{2}\left(  \left\langle x,\widehat{e}_{j}\right\rangle ^{2}-\left\langle x,e_{j}\right\rangle ^{2}\right)  }{\sum_{j=1}^{k_{n}}\lambda_{j}\left[  f_{n}\left(  \lambda_{j}\right)\right]  ^{2}\left\langle x,e_{j}\right\rangle ^{2}}.
\end{align*}
Copying what was done for the proof of Corollary \ref{Cor2}, we can easily prove that $w_{n1}\overset{\mathbb{P}}{\rightarrow}0.$ In order to alleviate formulas and displays, we are going to prove that $w_{n2}\overset{\mathbb{P}}{\rightarrow}0$ in the special case when $f_{n}\left(  \lambda_{j}\right)=1/\lambda_{j}.$ The general situation stems easily from this special case. Thus, we have now
$$
w_{n2}=\dfrac{\sum_{j=1}^{k_{n}}\left(  \left\langle x,\widehat{e}_{j}\right\rangle ^{2}-\left\langle x,e_{j}\right\rangle ^{2}\right)/\lambda_{j}}{\sum_{p=1}^{k_{n}}\left\langle x,e_{j}\right\rangle ^{2}/\lambda_{j}}.
$$
We denote by $\widehat{\pi}_{j}$ the projector on the eigenspace associated to the $j^{th}$ eigenfunction of $\Gamma_n$. Then, with this notation, we can write $\left\langle x,\widehat{e}_{j}\right\rangle ^{2}-\left\langle x,e_{j}\right\rangle ^{2}=\left\Vert\widehat{\pi}_{j}x\right\Vert _{\infty}^{2}-\left\Vert \pi_{p}x\right\Vert_{\infty}^{2}=\left\langle \left(  \widehat{\pi}_{j}-\pi_{j}\right)x;x\right\rangle $ and we have
$$
\left\vert \left\langle x,\widehat{e}_{j}\right\rangle ^{2}-\left\langle x,e_{j}\right\rangle ^{2}\right\vert \leq\left\Vert \widehat{\pi}_{j}-\pi_{j}\right\Vert _{\infty}\left\Vert x\right\Vert ^{2},
$$
\begin{align*}
\widehat{\pi}_{j}-\pi_{j}  &  \ =\ \dfrac{1}{2\pi\iota}\int_{\mathcal{B}_{j}}\left[  \left(  zI-\Gamma_{n}\right)  ^{-1}-\left(  zI-\Gamma\right)^{-1}\right]  dz\\
&  =\dfrac{1}{2\pi\iota}\int_{\mathcal{B}_{j}}\left[  \left(  zI-\Gamma_{n}\right)  ^{-1}\left(  \Gamma_{n}-\Gamma\right)  \left(  zI-\Gamma\right)^{-1}\right]  dz,
\end{align*}
and
$$
\mathbb{E}\left\Vert \widehat{\pi}_{j}-\pi_{j}\right\Vert _{\infty}\leq C\dfrac{j\log j}{\sqrt{n}}.
$$
Finally
$$
\left\vert w_{n2}\right\vert \leq C\dfrac{1}{\sqrt{n}}\sum_{j=1}^{k_{n}}j\log j\leq C\dfrac{k_{n}^{2}\log k_{n}}{\sqrt{n}}\rightarrow0,
$$
which finishes the proof of the Corollary.

\textbf{Proof of Proposition \ref{contrex2}}

Take $x=\sum x_{i}e_{i}$ and $\rho=\sum\rho_{i}e_{i}$ in $H$. 
Obviously, it suffices to prove that the Proposition holds when $\widehat{\Pi}_{k_{n}}-\Pi_{k_{n}}$ is replaced with
$$
\varphi_{k_{n}}\left(  \left(  \Gamma_{n}-\Gamma\right)  \right)  =\sum_{j=1}^{k_{n}}\left[  \mathcal{S}_{j}\left(  \Gamma_{n}-\Gamma\right)  \Pi_{j}+\Pi_{j}\left(  \Gamma_{n}-\Gamma\right)  \mathcal{S}_{j}\right].
$$
Following Dauxois {\it et al.} (1982), p. 143-144, we can check that when $X_{1}$ is Gaussian, $\sqrt{n}\left(  \Gamma_{n}-\Gamma\right)  $ converges weakly to the Gaussian random operator $G$ defined by
$$
G  =\sum_{j\leq j^{\prime}}\sqrt{\lambda_{j}\lambda_{j^{\prime}}}\xi_{j,j^{\prime}}\left(  e_{j}\otimes e_{j^{\prime}}+e_{j^{\prime}}\otimes e_{j}\right) +\sqrt{2}\sum_{j}\lambda_{j}\left(  e_{j}\otimes e_{j^{\prime}}\right)\xi_{j,j},
$$
where $\xi_{j,j^{\prime}}$'s are i.i.d. Gaussian centered r.r.v. with variance equal to 1. Thus, we replace once more $\sqrt{n}\left(  \Gamma_{n}-\Gamma\right)  $ with $G$ (the situation is indeed the same as if the operator $X_{1}\otimes X_{1}$ was assumed to be Gaussian). We are going to prove that $\dfrac{\left\langle \varphi_{k_{n}}\left(  G\right)  \rho,x\right\rangle }{t_{n,x}}$ is not bounded in probability, whatever the sequence $k_{n}\rightarrow+\infty,$ by choosing a special $\rho.$ We focus on the $jth$ term of the above sum.

The exact computation of $<\left(  \Pi_{j}G\mathcal{S}_{j}+\mathcal{S}_{j}G\Pi_{j}\right)  \left(  x\right)  ,\rho>$ may be deduced from Dauxois \textit{et al.} (1982) p.146. Assuming that all the $\lambda_{j}$'s have all the same order of multiplicity equals to $1$, we easily get
$$
 <\left(  \Pi_{j}G\mathcal{S}_{j}+\mathcal{S}_{j}G\Pi_{j}\right)  \left(x\right)  ,\rho> =\sum_{l\neq j}\frac{\sqrt{\lambda_{l}\lambda_{j}}}{\lambda_{j}-\lambda_{l}}\left(  x_{j}\rho_{l}+x_{l}\rho_{j}\right)  \xi_{jl}.
$$
The previous sum is a real centered Gaussian random variable with variance
$$
\sum_{l\neq j}\frac{\lambda_{l}\lambda_{j}}{\left(  \lambda_{j}-\lambda_{l}\right)  ^{2}}\left(  x_{j}\rho_{l}+x_{l}\rho_{j}\right)  ^{2}.
$$
Summing over $j$ provides the variance of $\left\langle \varphi_{k_{n}}\left(G\right)  \rho,x\right\rangle $
$$
\sum_{j=1}^{k_{n}}\sum_{l\neq j}\frac{\lambda_{l}\lambda_{j}}{\left(\lambda_{j}-\lambda_{l}\right)  ^{2}}\left(  x_{j}\rho_{l}+x_{l}\rho
_{j}\right)  ^{2}\geq\sum_{j=1}^{k_{n}}\lambda_{j}\rho_{j}^{2}\sum_{l\neq j}\frac{\lambda_{l}x_{l}^{2}}{\left(  \lambda_{j}-\lambda_{l}\right)  ^{2}}\geq\sum_{j=1}^{k_{n}}\lambda_{j}\rho_{j}^{2}\sum_{l=1}^{j-1}\frac{\lambda_{l}x_{l}^{2}}{\left(  \lambda_{j}-\lambda_{l}\right)  ^{2}}.
$$
For the sake of simplicity we assume that $x_{k}>0$ and $\rho_{k}>0$. Now if $x_{l}^{2}=l^{-1-\beta}$ and $\lambda_{l}=l^{-1-\alpha}$ the computation of the second sum stems from
$$
\sum_{l=1}^{j-1}\frac{\lambda_{l}x_{l}^{2}}{\left(  \lambda_{j}-\lambda
_{l}\right)  ^{2}}\sim{\displaystyle\int_{1}^{j-1}}\dfrac{s^{\alpha-\beta}}{\left(  1-\left(  \dfrac{s}{j}\right)  ^{1+\alpha
}\right)  ^{2}}ds\sim Cj^{2+\alpha-\beta}.
$$
Finally
$$
\sum_{j=1}^{k_{n}}\lambda_{j}\rho_{j}^{2}\sum_{l=1}^{j-1}\frac{\lambda_{l}x_{l}^{2}}{\left(  \lambda_{j}-\lambda_{l}\right)  ^{2}}\geq C\sum_{j=1}^{k_{n}}j^{1-\beta}\rho_{j}^{2}\rightarrow+\infty.
$$
We see that the variance of $\dfrac{\left\langle \varphi_{k_{n}}\left(G\right)  \rho,x\right\rangle }{t_{n,x}}$ explodes and that this random variable cannot converge in distribution.


\textbf{Proof of Theorem \ref{contrex}}


From (\ref{decomp}) and all that was made above it suffices to prove that the Theorem holds with $U_{n}$ replacing $\widehat{\rho}-\widehat{\Pi}_{k_{n}}\rho.$ Now suppose that for a given normalizing sequence $\alpha_{n}>0,$ $\alpha_{n}U_{n}$ converges weakly in the norm topology of $H$ (the deterministic sequence $\alpha_{n}$ just depends on the random variables $X_{k},\varepsilon_{k}$). For all $x$ in $H$, $\alpha_{n}\left\langle U_{n},x\right\rangle $ converges weakly too and
$$
\alpha_{n}\left\langle U_{n},x\right\rangle =\dfrac{\alpha_{n}}{n}\sum_{i=1}^{n}\left\langle X_{i},\Gamma^{\dag}x\right\rangle \varepsilon_{i},
$$
is an array of real independent random variable. Suppose that $x$ belongs to the domain of $\Gamma^{-1},$ namely that
$$
\sum_{j=1}^{+\infty}\dfrac{\left\langle x,e_{j}\right\rangle ^{2}}{\lambda_{j}^{2}}<+\infty,
$$
then
$$
\dfrac{1}{\sqrt{n}}\sum_{j=1}^{n}\left\langle X_{j},\Gamma^{\dag}x\right\rangle \varepsilon_{j}\overset{w}{\rightarrow}N\left(  0,\beta_{x}\sigma_{\varepsilon}^{2}\right)  ,
$$
where $\beta_{x}$ depends on $x$ and on the eigenvalues of $\Gamma.$ Consequently $\alpha_{n}=\sqrt{n}.$ Now if $\sum_{j}\left\langle x,e_{j}\right\rangle ^{2}/\lambda_{j}^{2}$ is divergent, $\mathbb{E}\left(  \left\langle X_{i},\Gamma^{\dag}x\right\rangle^{2}\varepsilon_{i}^{2}\right)  \uparrow+\infty$ and $\alpha_{n}\left\langle U_{n},x\right\rangle $ cannot converge in distribution. This finishes the proof of the Theorem.


\noindent{\bf Acknowledgements}. The authors would like to thank the organizers and participants of the working group STAPH on functional statistics  in Toulouse for fruitful discussions.

\end{document}